# THE HYPERBOLIC GEOMETRY OF RANDOM TRANSPOSITIONS[1]

By Nathanaël Berestycki

*University of British Columbia and Ecole Normale Supérieure, Paris*

Turn the set of permutations of $n$ objects into a graph $G_n$ by connecting two permutations that differ by one transposition, and let $\sigma_t$ be the simple random walk on this graph. In a previous paper, Berestycki and Durrett [In *Discrete Random Walks* (2005) 17–26] showed that the limiting behavior of the distance from the identity at time $cn/2$ has a phase transition at $c = 1$. Here we investigate some consequences of this result for the geometry of $G_n$. Our first result can be interpreted as a breakdown for the Gromov hyperbolicity of the graph as seen by the random walk, which occurs at a critical radius equal to $n/4$. Let $T$ be a triangle formed by the origin and two points sampled independently from the hitting distribution on the sphere of radius $an$ for a constant $0 < a < 1$. Then when $a < 1/4$, if the geodesics are suitably chosen, with high probability $T$ is $\delta$-thin for some $\delta > 0$, whereas it is always $O(n)$-thick when $a > 1/4$. We also show that the hitting distribution of the sphere of radius $an$ is asymptotically singular with respect to the uniform distribution. Finally, we prove that the critical behavior of this Gromov-like hyperbolicity constant persists if the two endpoints are sampled from the uniform measure on the sphere of radius $an$. However, in this case, the critical radius is $a = 1 - \log 2$.

**1. Introduction.** Let $\mathcal{S}_n$ be the set of permutations of $\{1, 2, \ldots, n\}$, and let $\sigma_t$ be the continuous-time random walk on $\mathcal{S}_n$ that results when randomly chosen transpositions are performed at rate 1. Let $d(\sigma_t)$ be the distance from the identity $I$ at time $t$, that is, the minimum number of transpositions needed to return to $I$. In a previous paper, Berestycki and Durrett [3] showed

Received December 2004; revised September 2005.

[1]Supported in part by Rick Durrett's joint NSF–NIGMS Grant DMS-02-01037.

THEOREM 0.  *As $n \to \infty$, $d(\sigma_{cn/2})/n \to_p u(c)$ where*

$$(1.1) \qquad u(c) = 1 - \sum_{k=1}^{\infty} \frac{1}{c} \frac{k^{k-2}}{k!} (ce^{-c})^k.$$

*Although it is not easy to see from the formula, the function $u(c) = c/2$ for $c \leq 1$ and is $< c/2$ for $c > 1$.*

We can think of $\sigma_t$ as a random walk on the graph $G_n$ with vertices $\mathcal{S}_n$ and edges connecting two permutations that differ by one transposition, so that $G_n$ is the Cayley graph of $\mathcal{S}_n$ associated with the set of generators $S = \{\text{all transpositions}\}$. Theorem 0 was proved by establishing a connection with Erdős–Rényi random graphs. The phase transition observed for $\sigma_t$ is then related to the well-known double jump of the size of connected components of $G(n, c/n)$ at $c = 1$. [Here and in all that follows, $G(n, p)$ denotes the Erdős–Rényi random graph with parameters $n$ and $p$, i.e., a random graph on $n$ vertices where each edge is present independently of the others with probability $p$.] We refer the reader to Janson et al. [8] for this and other facts about Erdős–Rényi random graphs.

In this paper we try to investigate some of the geometric implications of Theorem 0. We find a new connection between the speed of a random walk and the Gromov hyperbolicity of the space in which the random walk is evolving.

*Organization of the paper.*  In Sections 1.1, 1.2, 1.3 we present our results. The proofs of these results can be found successively in Sections 2–8. Each proof is preceded by a restatement of the corresponding theorem for convenience, and by an informal proof which outlines the main ideas used.

1.1. *Asymptotic hyperbolicity.*  The notion of hyperbolicity for a discrete structure such as a group is a notion that goes back to Gromov [7]. As there is no derivative, and thus no curvature available in a discrete space, the idea is to define what hyperbolic means using only elementary properties of the space.

One way to do this is as follows. Let $(X, |\cdot|)$ be a metric space, where $|x - y|$ denotes the distance between $x$ and $y$. For points $x, y$ and $p$ in $X$, define the Gromov inner product by

$$2(x|y)_p = |x - p| + |y - p| - |x - y|.$$

$(x|y)_p$ thus measures how well the union of the geodesic segments $[p, x] \cup [p, y]$ approximates a geodesic between $x$ and $y$. Gromov's original definition of hyperbolic spaces is as follows. Call $X$ $\delta$-hyperbolic if

$$(1.2) \qquad (x|z)_p \geq (x|y)_p \wedge (y|z)_p - \delta$$

for all $x, y, z$ and $p$. This definition is not very intuitive at first, but fortunately there is an equivalent definition, which can be formulated using the notion of $\delta$-thin triangle. A triangle $(x, y, z)$ with geodesic sides $s_1, s_2, s_3$ is said to be $\delta$-thin if any side, say $s_1$, lies entirely within distance at most $\delta$ of the two remaining sides:

$$s_1 \subset \{x \in X, d(x, s_2 \cup s_3) \le \delta\}.$$

The space is called $\delta$-hyperbolic if all geodesic triangles are $\delta$-thin, and it is simply called hyperbolic if it is $\delta$-hyperbolic for some $\delta \ge 0$ (when $\delta = 0$, the space isometrically embeds into a tree). It is not immediate, but not hard to check, that if *all* triangles $(x, y, z)$ are $\delta$-thin, then (1.2) is satisfied for some number $\delta'$ that may differ by a constant factor from $\delta$. Conversely, in a space where (1.2) is satisfied for *all* points $(p, x, y, z)$, all triangles are $\delta'$-thin, where $\delta'$ may differ from $\delta$ by a constant factor.

Of course a bounded space (in particular, a finite space such as $\mathcal{S}_n$) is trivially hyperbolic, but we will be interested in situations where the constant $\delta$ may or may not stay bounded as the size of the space tends to $\infty$.

Our first result makes the connection between Theorem 0 and Gromov hyperbolic spaces, where we look at the two definitions of hyperbolic constants suitably weakened. For $0 < a < 1$, let $\partial B(an)$ be the sphere of radius $an$, that is, the set of points at distance $\lfloor an \rfloor$ from the origin. We let $\nu$ be the hitting distribution of $\partial B(an)$ by $\sigma_t$, that is, $\nu$ is the law on $\partial B(an)$ of $\sigma_T$ where $T = \inf\{t > 0, d(\sigma_t) = \lfloor an \rfloor\}$.

THEOREM 1. *Let $x, y$ be sampled from $\nu$ independently, and set $p = I$, the identity element.*

1. *If $a < 1/4$, then there is some $\delta < \infty$ (depending only on $a$), such that*

$$E(x|y)_p \le \delta.$$

   *Moreover, with probability asymptotically 1, there is a geodesic between $x$ and $y$ that comes within expected distance $\delta' < \infty$ of $p$.*

2. *If $a > 1/4$, then*

$$E(x|y)_p \sim \delta n$$

   *for some $0 < \delta < \infty$. Moreover, no geodesic between $x$ and $y$ can approach $p$ closer than $\delta'n$ with probability asymptotically 1, where $0 < \delta' < \infty$.*

In the statement of the theorem and in the rest of the paper, $a_n \sim b_n$ means that $a_n/b_n \to 1$.

REMARK. It follows immediately from Theorem 1 that when $a < 1/4$, with probability asymptotically 1,

$$(x|z)_p \ge (x|y)_p \wedge (y|z)_p - \delta$$

for independent $x, y, z$ sampled from $\nu$, hence the idea that definition 1 of hyperbolicity is satisfied "asymptotically $\nu$-almost surely." The statement about the geodesics shows that definition 2 is satisfied "asymptotically $\nu$-almost surely" when $a < 1/4$.

At this point we should emphasize that the result in Theorem 1 involves hyperbolic constants that are different from the standard definitions discussed above in several important ways. The most obvious difference comes from the randomness of $x$ and $y$, and from the fact that the roles played by $x, y$ and $p$ are somewhat different. Here $p$ is a fixed reference point, whereas Gromov's definition requires that every triangle should be thin. Another issue is that, corresponding to the second definition of hyperbolicity with thin triangles, we show that there exists a certain geodesic between $x$ and $y$ having the desired properties. As we will see below in Theorem 6, there may be a great many geodesics between two given points in $\mathcal{S}_n$. More importantly, these geodesics can be far apart, as will show the following concrete example:

$$
\begin{array}{ll}
\sigma & (1\ \ 14\ \ 5\ \ 11)\ \ (2)\ \ (3\ \ 9)\ \ (4\ \ 13\ \ 6)\ \ (7\ \ 12\ \ 8\ )\ \ (10) \\
\pi_1 & (1)\ \ (14)\ \ (5)\ \ (11)\ \ (2)\ \ (3)\ \ (9)\ \ (4\ 13\ 6)\ \ (7\ \ 12\ \ 8\ )\ \ (10) \\
\pi_2 & (1\ \ 14\ 5\ \ 11)\ \ (2)\ \ (3\ \ 9)\ \ (4)\ \ (13)\ \ (6)\ \ (7)\ \ (12)\ \ (8)\ \ (10) \\
\pi_1\pi_2^{-1} & (11\ \ 5\ \ 14\ \ 1)\ \ (2)\ \ (9\ \ 3)\ \ (4\ \ 13\ \ 6)\ \ (7\ \ 12\ \ 8)\ \ (10)
\end{array}
$$

Since for any permutation $\pi$ we have $d(\pi) = n - \#$ cycles of $\pi$, $d(\sigma) = 8$. $\pi_1$ and $\pi_2$ are on two geodesics from $I$ to $\sigma$, but $d(\pi_1, \pi_2) = d(\pi_1\pi_2^{-1}) = 8$. In general if $d(\sigma) = cn/2$ with $c < 1$, and we divide the cycles at random into two groups, we can define $\pi_1$ to have cycle structure given by the first group of $\sigma$ staying as it is and the second completely broken in cycles on lengths 1. If we define $\pi_2$ by the exchange of the two groups, then we will have $d(\sigma, \pi_i) = cn/4$ and $d(\pi_1, \pi_2) = cn/2$.

1.2. *The geometry of* $G_n$.   How much can we learn from Theorem 1 about the global geometry of $G_n$? To answer this question, we need to see how special a choice it is to sample the points $x$ and $y$ according to the hitting distribution $\nu$. (The fact that $p = I$ is a fixed reference point is not too important, due to the transitivity of $G_n$.) We begin by an apparently unrelated question, which is to ask how large is a ball of radius $an$.

THEOREM 2.   *If* $0 \le a \le 1$, *then as* $n \to \infty$, *we have* $|B(I, an)| \approx (n!)^a$ *in a logarithmic sense, that is,*

$$
\lim_{n \to \infty} \frac{\log |B(I, an)|}{n \log n} = a.
$$

This result is probably not new, but we have not found it in the literature. Our original motivation for studying the volume growth in $G_n$ was to try to understand the phase transition of Theorem 0 in terms of the geometry of $G_n$. Our first thought was that since the speed was nonsmooth we might see a change in the volume growth. The above result contradicts this idea.

To put our next two results into perspective it is useful to contrast them with Brownian motion $B_t$ on a $d$-dimensional manifold of constant negative curvature $-1$. In that case as $t \to \infty$, if $d(B_t)$ is the distance from the origin, then (see [11], e.g.) there is a constant $v$ so that

$$d(B_t)/t \to v \qquad \text{as } t \to \infty.$$

In the case of Brownian motion on hyperbolic space, rotational symmetry implies that the hitting distribution is uniform. In contrast for the random transposition random walk, we will see in Theorem 3 that the hitting distribution is asymptotically singular with respect to the uniform distribution on $\partial B(I, an)$.

THEOREM 3.    *Let $|\mathcal{C}_1|$ be the length of the cycle that contains $1$. Under $\mu$, the uniform distribution on $\partial B(I, an)$,*

$$|\mathcal{C}_1| \Rightarrow \mathbf{G}$$

*where $\mathbf{G}$ is a geometric r.v. with $P(\mathbf{G} > k) = (b/(1+b))^k$ and $b$ satisfies $\log(1+b)/b = 1 - a$.*

To describe the hitting distribution $\nu$, we note that (1.1) suggests that it will be the same as the distribution of $\sigma_{cn/2}$ where $c = u^{-1}(a)$. When $a > 1/2$ this is much different from the distribution in Theorem 3 since in this case $c > 1$ and Schramm [12] has shown that $\sigma_{cn/2}$ has cycles of lengths of order $n$.

Here we will concentrate on what happens when $a < 1/2$ and $c = 2a$. In this case results in [3] show that as $n \to \infty$, the number of fragmentations before time $cn/2$ is asymptotically a Poisson random variable with mean $\kappa(c) = -(\log(1-c) + c)/2$. In particular,

$$P(d(\sigma_{cn/2}) = cn/2) \to e^{-\kappa(c)} = e^{c/2}\sqrt{1-c}.$$

It will be convenient to approach the hitting distribution $\nu$ by the distribution $\nu_0$ of $\sigma_{cn/2}$ conditioned on no fragmentation. More generally, if $\nu_k = \nu$ conditioned on exactly $k$ fragmentations before the hitting time,

$$\nu = e^{-\kappa(c)} \sum_{k=0}^{\infty} \nu_k \frac{\kappa(c)^k}{k!} + o(1).$$

To study $\nu_0$, we recall the connection with random graphs developed in [3]: when we transpose $i$ and $j$ we draw an edge between $i$ and $j$. In order for the distance from the identity to increase by 1 at each time, each transposition must involve indices from two different cycles and will merge them into one. In terms of the random graph, this means that all components are trees. Using results from [3], it is straightforward to show:

THEOREM 4. *Let $\mathcal{C}_1$ be the length of the cycle that contains 1. Let $c < 1$. Under $\nu_0$,*

$$P(|\mathcal{C}_1| = k) \to \frac{1}{c} \frac{k^{k-1}}{k!} (ce^{-c})^k \qquad \text{for all } k \geq 1.$$

Theorems 3 and 4 show that the uniform distribution $\mu$ and the hitting distribution $\nu_0$ concentrate on different permutations. In the first case the number of fixed points will be close to its expected value $nP(|\mathcal{C}_1| = 1) = n/(1+b)$. In the second it will be close to $ne^{-c}$ by Theorem 4. This is made precise by the following theorem.

THEOREM 5. *As $n \to \infty$, the hitting distribution $\nu$ and the uniform distribution $\mu$ on a sphere of radius $an$ are asymptotically singular:*

$$d_{\mathrm{TV}}(\mu, \nu) \to 1.$$

Let $t = [cn/2]$ with $c < 1$. To understand why $\nu$ is different from $\mu$ we will examine the Radon–Nikodym derivative $r(\sigma) = d\nu_0/d\mu$. It is not hard to show that

THEOREM 6. *Suppose $d(\sigma) = t$ and $m_1, \ldots, m_j$ are the cycle lengths of $\sigma$. The number of paths of length $t$ from $I$ to $\sigma$ is*

$$t! \prod_{i=1}^{j} \frac{m_i^{m_i-2}}{(m_i-1)!}.$$

*If $t = \lfloor cn/2 \rfloor$ with $c < 1$, then*

$$r(\sigma) = K_{n,t} \prod_{i=1}^{j} \frac{m_i^{m_i-2}}{(m_i-1)!},$$

*where $K_{n,t}$ is a constant that only depends on $n$ and $t$.*

The last result enables us to prove a stronger version of Theorem 5: it tells us that the "support" of $\nu$ is concentrated on a set that is exponentially smaller than the size of $\partial B(an)$.

THEOREM 7. *Suppose $a < 1/2$. There exists a set $S_n \in \partial B(an)$ such that $\nu(S_n) \to 1$ as $n \to \infty$ and*

$$\lim_{n \to \infty} \frac{1}{n} \log \frac{|S_n|}{|\partial B(an)|} = \gamma < 0.$$

1.3. *The hyperbolic constant under the uniform measure.* In Theorem 1, we learn that if $x$ and $y$ are sampled from $\nu$, roughly speaking, the Gromov hyperbolicity of the "support" breaks down at $a = 1/4$, that is, the hyperbolic constant increases suddenly from $O(1)$ to $O(n)$ at this point.

However, the results from the previous section tell us that this "support" is (exponentially) small with respect to the ambient space. It is therefore natural to ask what happens to Theorem 1 when we replace $\nu$ with the uniform measure $\mu$ on $\partial B(an)$. Theorem 8 will show that the qualitative behavior of the hyperbolic constant remains the same. We prove that there is a threshold where the expected Gromov inner product $E(\sigma|\pi)_p$ jumps from $O(1)$ to $O(n)$, but this time the critical value is $a = 1 - \log 2 \approx 0.31$, rather than $a = 1/4$.

When $\sigma$ and $\pi$ are independent uniform permutations on $\partial B(an)$, by the transitivity of $G_n$, it is enough to analyze $d(\sigma, \pi)$ to understand $(\sigma|\pi)_p$, the inner Gromov product. Since $d(\sigma, \pi) = d(I, \sigma^{-1}\pi)$, which has the same law as $d(I, \sigma\pi)$, it will be enough to characterize the values of $a$ for which $d(I, \sigma\pi) = 2an + o(n)$ and those for which it is $< 2an$.

THEOREM 8. *Let $0 < a < 1$ and let $\sigma, \pi$ be two random independent points chosen uniformly from $\partial B(an)$. Then:*

1. *If $a < 1 - \log 2$,*

$$E(\sigma|\pi)_p \leq \delta(\log n)^2$$

   *for some $0 < \delta = \delta(a) < \infty$. Moreover, with probability asymptotically 1, there is a geodesic between $\sigma$ and $\pi$ that comes within distance at most $\delta(\log n)^2$ of $p$.*
2. *If $a > 1 - \log 2$,*

$$E(\sigma|\pi)_p \sim \delta n$$

   *for some $\delta = \delta(a) > 0$. Moreover, no geodesic can approach $p$ closer than $\delta' n$ for some $0 < \delta' < \infty$.*

REMARK. The $O((\log n)^2)$ bound in part 1 of the theorem could probably be improved into an $O(1)$ bound (just like in Theorem 1) with some more work, but we have not tried to do so. In part 2, by analogy with Berestycki and Durrett [3], we conjecture that the fluctuations are of order exactly $n^{1/2}$ in the supercritical regime. More precisely, it should be true that when $a > 1 - \log 2$,

$$n^{-1/2}(E(\sigma|\pi)_p - \delta n) \Rightarrow \mathcal{N}(0, \kappa),$$

where $\delta$ is the limit in part 2 of the theorem, and $\kappa$ is some parameter.

**2. Asymptotic hyperbolicity under $\nu$.** The first result we prove is Theorem 1.

THEOREM 1. *Let $x, y$ be sampled from $\nu$ independently, and set $p = I$, the identity element.*

1. *If $a < 1/4$, then there is some $\delta < \infty$ (depending only on $a$), such that*

$$E(x|y)_p \leq \delta.$$

   *Moreover, with probability asymptotically 1, there is a geodesic between $x$ and $y$ that comes within expected distance $\delta' < \infty$ of $p$.*

2. *If $a > 1/4$, then*

$$E(x|y)_p \sim \delta n$$

   *for some $0 < \delta < \infty$. Moreover, no geodesic between $x$ and $y$ can approach $p$ closer than $\delta' n$ with probability asymptotically 1, where $0 < \delta' < \infty$.*

*Sketch of the proof.* Let $X_t$ and $Y_t$ be two independent random walks starting at the origin. Let them run until the times $T$ and $T'$ where they respectively hit the sphere $\partial B(an)$. Then the transitivity of the Cayley graph of $\mathcal{S}_n$, and the reversibility of the increments of the random walk, imply that $(X_T, X_{T-1}, \ldots, p, Y_1, \ldots, Y_{T'})$ is a random walk path of length $T + T'$. Hence the distance between $X_T = x$ and $Y_{T'} = y$ is the same as $d(\sigma_{T+T'})$. By Theorem 0, $T$ and $T' \approx \frac{1}{2}u^{-1}(a)n$, so applying Theorem 0 again, when $a < 1/4$, $|x - y| \approx 2an = |x| + |y|$ [the random walk runs for a time $2an < n/2$ and there are only $O(1)$ fragmentations]. For $a > 1/4$, the random walk is run for time $u^{-1}(a)n$ which, in view of Theorem 0, means that $c = 2u^{-1}(a)$, and $|x - y| = nu(2u^{-1}(a)) \ll 2an$. See Figure 1.

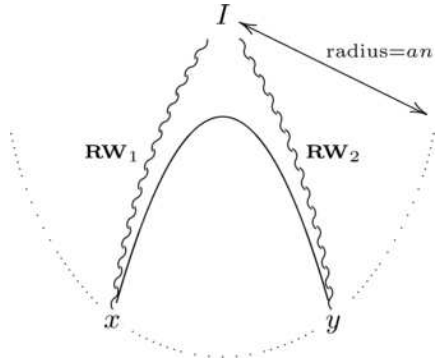

FIG. 1. *Two independent random walks run until they hit the sphere of radius $an$.*

The claim about the existence of a geodesic that makes the triangle $(x, p, y)$ thin involves necessarily another argument, since geodesics may be far apart. However, it is not very hard to construct by hand a geodesic between the identity and $x$ such that each point of the random walk path is within $O(1)$ of this geodesic. Applying this construction to the two random walk paths gives the result of Theorem 1.

PROOF OF THEOREM 1. Let us first deal with case $a < 1/4$ and prove that in this case $E(x|y)_p \leq \delta$. Keeping the same notation as above, note that $an + O(1)$ steps are sufficient for $X$ to reach distance $an$. Indeed, after $an$ steps, $X_{an}$ is at distance $an - X_1$ where $X_1$ is twice a Poisson random variable by Theorem 1 in [3]. It is immediate that in $X_1$ steps the probability that $X$ has a fragmentation converges to 0. Therefore $T - an \Rightarrow X_1$ (remember that here time is measured discretely). Similarly $T' - an \Rightarrow X_2$ where $X_1, X_2$ are i.i.d. Hence $(x = X_t, X_{t-1}, \ldots, X_1, I, Y_1, \ldots, Y_t = y)$ is a random walk of $2an + X_1 + X_2$ steps. In the worst case possible all $X_1 + X_2$ steps represent "backward" steps (meaning, toward $x$ rather than $y$). Hence if $X_3 = an - d(I, \sigma_{2an})$ (so that $X_3$ is also twice a Poisson random variable, but with a different parameter),

$$2E(x|y)_p = 2an - |x - y|$$
$$\leq 2an - Ed(I, \sigma_{2an}) + E(X_1 + X_2)$$
$$\leq E(X_1) + E(X_2) + E(X_3) < \infty.$$

It is slightly simpler to prove that when $a > 1/4$, $E(x|y)_p \sim \delta n$. Indeed, in this case, by Theorem 0, we have that

$$\tfrac{1}{2} u^{-1}(a) - \varepsilon \leq T/n \leq \tfrac{1}{2} u^{-1}(a) + \varepsilon.$$

Therefore

$$\inf_{|t/n - u^{-1}(a)| \leq 2\varepsilon} d(I, \sigma_t) \leq |x - y| \leq \sup_{|t/n - u^{-1}(a)| \leq 2\varepsilon} d(I, \sigma_t).$$

An easy estimate shows that we are never off by more than $O(n^{1/2})$ if we evaluate the distance of the random walk by counting the number of clusters of the random graph rather than the number of cycles of $\sigma_t$. But for the random graph, the number of clusters is monotone increasing. Hence, if $\alpha$ denotes $u(2u^{-1}(a))$, we have by continuity of $u$ that

$$\alpha - \varepsilon' + o(1) \leq \frac{E|x - y|}{n} \leq \alpha + \varepsilon' + o(1)$$

and $\varepsilon'$ can be made as small as desired by continuity of $u$. Therefore

$$\frac{E|x - y|}{n} \to \alpha.$$

It suffices now to prove that $\alpha < 2a$, that is, $u(2u^{-1}(a)) < 2a$ or, after change of variable $c = u^{-1}(a)$, it suffices to prove $u(2c) < 2u(c)$ for all $c > 1/2$. This fact is a consequence of the sublinearity of $u$: it will be proved later that $u$ is strictly concave on $[1, \infty)$, from which it follows that $u(c) > u(2c)/2$.

We now turn to the part of the theorem that concerns geodesics, and prove that for a random walk $(X_t, t \leq cn/2)$ of time-duration $cn/2$ with $c < 1$, there is a geodesic between $\sigma = X_{cn/2}$ and $I$, that we call $\gamma$, such that

$$(2.1) \qquad E \sup_{t \leq cn/2} d(X_t, \gamma) = O(1).$$

This shows that when $c < 1$ there is a geodesic that stays close to the random walk path. When $a < 1/4$, $p = I$ is on the random walk path that leads from $x$ to $y$, so this shows that $E(d(I, \gamma)) = O(1)$, as claimed in the theorem. The case $a > 1/4$ is trivial by the triangle inequality.

Let $\tau_1, \ldots, \tau_N$ be the sequence of transpositions that are the increments of the random walk path leading to $\sigma$, so that $\sigma = \tau_1 \ldots \tau_N$. Let $\gamma$ be the geodesic between $\sigma$ and $I$ defined by $\gamma_0 = \sigma$, $\gamma_1 = \sigma\tau_N$, $\gamma_2 = \gamma_1\tau_{N-1}, \ldots$, until the first time $t$ such that multiplying $\gamma_t$ by $\tau_{N-t}$ would result in a coagulation of two cycles of $\gamma_t$. We do not allow this possibility (otherwise $\gamma$ would not be a geodesic), and simply skip $\tau_{N-t}$: $\gamma_{t+1} = \gamma_t\tau_{N-t-1}$. We will see in a moment that this path never backtracks and that it ends at a bounded distance from $I$, to which it will be necessary to add a (bounded) number of steps so that it actually ends at $I$.

Let $n(t)$ be the index of the transposition to be performed at time $t$ on $\gamma_t$. Note that we can always write

$$\gamma_t = \tau_1\tau_2 \ldots \tau_{n(t)} \prod_{i \in K_t} \tau_i,$$

where $K_t$ is a set whose size we will show is bounded. Indeed, even when we skip $\tau_{n(t)}$ in $\gamma_t$ [so that $n(t) \in K_{t+1}$], the following transpositions $\tau_{n(t)-1}, \ldots$ commute with the members of $K_t$ with high probability and they can "jump above" the terms in $K_t$ and cancel the rest of the transpositions $(\tau_1 \ldots \tau_{n(t)-1})$.

LEMMA 1. *For all $t$, $E(|K_t|) \leq O(1)$, where $O(1)$ is a constant that depends only on $c < 1$. As a consequence, the path ends at bounded distance from the identity and the distance $E\sup_{t \leq cn/2} d(X_t, \gamma) = O(1)$.*

PROOF. There are two ways to add a member to $K_{t-1}$ at time $t$. The first one is that performing $\tau_{n(t)}$ will result in a coagulation, so that it is skipped by $\gamma$. The other way is if $\tau_{n(t)}$ does not commute with one of the members of $K_{t-1}$, it stays stuck somewhere in $K_t$.

If $\tau_{n(t)} = (i, j)$, we claim that in order for $i$ and $j$ to be in the same cycle of $\gamma_t$, $i$ and $j$ must belong to a component of the Erdős–Renyi graph associated with the random walk that contains a cycle at time $cn/2$. We will prove this in a moment, but if we admit this, then it follows that all transpositions in $K_t$ act on vertices that belong to $U(cn/2)$, the unicyclic components of the random graph at time $cn/2$: if $i \in K_t$, then either $\tau_i = (i,j)$ yields a coagulation in $\gamma_t$, or it does not commute with some member of $(k,l)$ of $K_{t-1}$, in which case $(i,j)$ overlaps with $(k,l)$. By induction, $k,l \in U(cn/2)$, therefore so are $i$ and $j$.

Let us prove our claim that if $(i,j)$ would yield a coagulation in $\gamma_t$, then $i,j \in U(cn/2)$. Let us observe first that $i$ and $j$ must already be in the same component of the random graph: because $\tau_{n(t)}$ was performed on the random walk, $i$ and $j$ were connected at that point in the random graph and they remain so. If $i$ and $j$ are in different cycles of $\sigma$, then there must have been some ulterior fragmentation in their cycles, so the claim holds. When they are in the same cycle of $\sigma$, then there must be some transposition $\tau_i$ with $i \in K_t$ such that $i$ and $j$ are in different cycles of $\gamma$ after $\tau_i$. Call those cycles $C_1$ and $C_2$. $\tau_i$ involves two members $k$ and $l$ of $C_1 \cup C_2$. Moreover the cycle structure of $\gamma$ before $(k,l)$ is performed must be of the form

$$(k, \ldots, i, \ldots, j, \ldots, l, \ldots);$$

otherwise $(k,l)$ cannot separate $i$ and $j$ at the next step. Unless $i$ and $j$ belong to a complex component, this implies that the cycle structure of $\sigma$ has the same form. However, this can only happen if $k$ and $l$ were connected to the component of $i$ and $j$ at different times; otherwise the cycle structure would be of the form $(i, \ldots j, \ldots, k, \ldots, l)$ or $(i, \ldots k, \ldots, l, \ldots, j)$. This implies in turn the existence of a cycle in the random graph component of $i$ and $j$ at time $cn/2$.

From there it follows in a straightforward way that $|K_t| \leq |U(cn/2)|$ (in a unicyclic component there are as many edges as vertices). It is now standard in the theory of random graphs to show that $|U(cn/2)|$ is bounded:

$$E|U(cn/2)| \sim \left(\frac{\pi}{8}\right)^{1/2} \sum_{k=2}^{\infty} \binom{n}{k} \frac{k^{k-1/2}}{k!} \left(\frac{c}{n}\right)^k \left(1 - \frac{c}{n}\right)^{k(n-k)+\binom{k}{2}-k} k$$

$$\sim \left(\frac{\pi}{8}\right)^{1/2} \sum_{k=2}^{\infty} \frac{k^{k+1/2}}{k!} (ce^{-c})^k < \infty$$

which completes the proof of the lemma. $\square$

Now let $X_t$ be a point on the random walk path. Since $\gamma$ tries to perform all $\tau_i$ (at reverse), there is a time $s$ such that $n(s) = t$, that is, the next transposition to be examined by $\gamma_s$ is $\tau_t$. At this time,

$$\gamma_s = \tau_1 \cdots \tau_t \prod_{i \in K_s} \tau_i$$

so that $|K_s|$ steps are enough to reach $\gamma_s$ from $X_t$. Since $E(|K_s|) < O(1)$ by the lemma, we have proved that

$$E \sup_{t \leq cn/2} d(X_t, \gamma) \leq O(1)$$

and Theorem 1 is proved. □

**3. Large deviations and volume growth.** The goal of this section is to prove Theorem 2, which we restate here for convenience.

THEOREM 2. *If* $0 \leq a \leq 1$, *then as* $n \to \infty$, *we have* $|B(I, an)| \approx (n!)^a$ *in a logarithmic sense, that is,*

$$\lim_{n \to \infty} \frac{\log |B(I, an)|}{n \log n} = a.$$

*Sketch of the proof.* The proof of the result is more interesting than the limit. We begin by recalling the dynamics of the Chinese restaurant process (see, e.g., [9]). Customer 1 enters and sits at table 1. At step $i$, customer $i$ enters and starts a new table with probability $1/i$ or sits to the left of customer $k$ where $k$ is chosen uniformly at random in $\{1, \ldots, i\}$. From the tables we define a permutation $\sigma$ by $\sigma(i) = i$ if customer $i$ is sitting by himself at his table and $\sigma(i) = k$ if $k$ sits to the right of $i$. It is easy to see that this defines a uniform random permutation on $\mathcal{S}_n$, and that the cycle structure is given by listing the individuals at the tables in clockwise order. It is well known that if $\sigma \in \mathcal{S}_n$, then $d(\sigma) = n-$ the number of cycles of $\sigma$. In the Chinese restaurant process construction, let $\zeta_i$ be the random variables taking the value 1 if customer $i$ sits at an existing table (and 0 otherwise). The $\zeta_i$'s are independent Bernoulli$(1 - 1/i)$ random variables. Recall that if $\sigma$ is a permutation, then $d(\sigma, I) = d(\sigma)$ is $n - \#$cycles of $\sigma$. Hence, if $\sigma$ is uniformly distributed over $\mathcal{S}_n$, then $d(\sigma)$ has the same distribution as $S_n = \sum_{i=1}^n \zeta_i$.

The $i$'s where a new cycle starts (i.e., $\zeta_i = 0$) are distributed with the same law as that of the occurrences of records for i.i.d. variables with continuous distribution function (cf. [5], Example 6.2 of Chapter 1). From calculations in that example it follows that $(n - S_n)/\log n \to 1$ in probability.

Returning to our calculation of the volume of the ball,

$$|B(I, an)| = n! P(S_n \leq an)$$

for all $0 < a < 1$. It is straightforward to generalize large deviation results for i.i.d. random variables (see, e.g., [5], Section 2.9) to prove Theorem 1. One begins with the observation that for $\lambda > 0$

$$(3.1) \qquad P(S_n \leq an) \leq e^{\lambda an} E e^{-\lambda S_n}$$

optimizes the upper bound over $\lambda$ and uses a change of measure argument to prove a corresponding lower bound.

PROOF OF THEOREM 2. Let $\{\zeta_i, \ i \geq 1\}$ be independent with $P(\zeta_i = 1) = 1 - 1/i$, and let $S_n = \sum_{i=1}^{n} \zeta_i$. Since $(\log n!)/(n \log n) \to 1$ it suffices to show:

LEMMA 2. Let $0 < a < 1$. As $n \to \infty$,

$$\lim_{n \to \infty} \frac{\log P[S_n \leq an]}{n \log n} = a - 1.$$

PROOF. Let $\varphi_n(\lambda) = \mathbf{E}[e^{-\lambda S_n}]$. Using the definition we have

$$\varphi_n(\lambda) = \prod_{i=1}^{n} \left[ \left( 1 - \frac{1}{i} \right) e^{-\lambda} + \frac{1}{i} \right] \equiv \prod_{i=1}^{n} q_i,$$

where $\equiv$ indicates that the last equation is the definition of $q_i$. By Markov's inequality we have

$$(3.2) \qquad \log P[S_n \leq an] \leq n \left( \lambda a + \frac{1}{n} \log \varphi_n(\lambda) \right) \qquad \text{for all } \lambda.$$

If we define

$$F_\lambda(x) = \frac{1}{\varphi_n(\lambda)} \int_{-\infty}^{x} e^{-\lambda y} \, dF^n(y),$$

then $F_\lambda$ is a distribution function such that

$$\text{mean}(F_\lambda) = -\frac{\varphi'_n(\lambda)}{\varphi_n(\lambda)} \quad \text{and} \quad \text{var}(F_\lambda) = \frac{d}{d\lambda} \frac{\varphi'_n(\lambda)}{\varphi_n(\lambda)} \geq 0.$$

To optimize (3.2), we want to choose $\lambda$ so that

$$a + \frac{1}{n} \frac{\varphi'_n(\lambda)}{\varphi_n(\lambda)} = 0.$$

This says that the mean of the transformed distributions is $na$, so

$$a = \frac{1}{n} \sum_{i=1}^{n} \frac{(1 - 1/i)e^{-\lambda}}{q_i} = \frac{1}{n} \sum_{i=1}^{n} \frac{(i-1)e^{-\lambda}}{(i-1)e^{-\lambda} + 1}.$$

We guess that the optimal $\lambda$ must be given (asymptotically) by $e^{-\lambda_{\text{opt}}} = b/n$. Plugging this in the above gives

$$a = \frac{1}{n} \sum_{j=1}^{n-1} \frac{jb/n}{(jb/n) + 1} \to \int_{0}^{1} \frac{bx}{bx + 1} \, dx = 1 - \frac{1}{b} \log(b + 1).$$

From this we see that we should choose $b$ so that $\log(b + 1)/b = 1 - a$.

*Upper bound.* Let us calculate what (3.2) gives with this choice of $\lambda$:

$$(3.3) \qquad \frac{1}{n}\log\varphi_n(\lambda) = -\log n + \frac{1}{n}\log\prod_{i=1}^{n}\left(\left(1-\frac{1}{i}\right)b + \frac{n}{i}\right)$$

$$(3.4) \qquad \qquad \to -\log n + \int_0^1 \log(b+1/x)\,dx.$$

Since the last integral is finite it follows from (3.2) and $\lambda_{\mathrm{opt}} = -\log b + \log n$ that

$$\limsup_{n\to\infty}\frac{1}{n\log n}P(S_n \le na) \le a - 1,$$

proving the upper bound half of Lemma 2.

*Lower bound.* The argument is similar to that in [5], page 73. Fix any $\nu < a$ and $\nu < \nu' < a$. Define a real number $b'$ by

$$\frac{\log(1+b')}{b'} = 1 - \nu'.$$

For any $\lambda$,

$$P[S_n \le an] \ge \int_{\nu n}^{an} dF^n(x) \ge \int_{\nu n}^{an} e^{\lambda x}\varphi_n(\lambda)\,dF_\lambda(x)$$

$$\ge \varphi_n(\lambda)e^{\lambda n\nu}[F_\lambda(na) - F_\lambda(n\nu)].$$

First, we prove that we can choose $\lambda$ such that $[F_\lambda(na) - F_\lambda(n\nu)] \to 1$. Recall that the mean of $F_\lambda$ is $-\frac{\varphi'_n(\lambda)}{\varphi_n(\lambda)}$, and that the latter function starts at $n - \log n$ for $\lambda = 0$, is strictly decreasing and equals $na$ when $\lambda = \lambda_{\mathrm{opt}} = -\log b + \log n$, that is, $e^{-\lambda_{\mathrm{opt}}} = b/n$. Thus if we pick $\lambda = \lambda'$ such that $e^{-\lambda'} = b'/n$, the mean of $F_{\lambda'}$ is by the lower bound calculation exactly $n\nu'$, and we have chosen $\nu < \nu' < a$. To conclude that $F_{\lambda'}(na) - F_{\lambda'}(n\nu) \to 1$, instead of using a law of large numbers argument such as in the i.i.d. case, we simply compute the variance of $F_{\lambda'}$ directly. Anticipating on the calculations of the next section, breaking the factor $e^{-\lambda x}$ in the Radon–Nikodym derivative of $F_\lambda$ into $e^{-\lambda\sum x_i}$ means that we can see $F_\lambda$ as a sum of independent Bernoulli random variables with parameter $\beta_i$ so that the variance is

$$\mathrm{var}\,F_{\lambda'} = \sum_{i=1}^{n}\beta_i(1-\beta_i) \le \sum_{i=1}^{n}\beta_i = \mathrm{mean}(F_\lambda) = n\nu' = O(n).$$

Another way to obtain this inequality is to do more direct computations:

$$\mathrm{var}\,F_\lambda = \frac{\varphi''_n(\lambda)}{\varphi_n(\lambda)} - \left(\frac{\varphi'_n(\lambda)}{\varphi_n(\lambda)}\right)^2$$

$$= \sum_{i=1}^{n} \frac{e^{-\lambda}(1-1/i)}{q_i} + \sum_{i \neq j} \frac{e^{-\lambda}(1-1/i)e^{-\lambda}(1-1/j)}{q_i q_j}$$

$$- \left( \sum_{i=1}^{n} \frac{e^{-\lambda}(1-1/i)}{q_i} \right)^2$$

$$= \sum_{i=1}^{n} \frac{e^{-\lambda}(1-1/i)}{q_i} - \sum_{i=1}^{n} \frac{e^{-2\lambda}(1-1/i)^2}{q_i^2}$$

$$\leq \sum_{i=1}^{n} \frac{e^{-\lambda}(1-1/i)}{q_i} = \frac{\varphi_n'(\lambda)}{\varphi_n(\lambda)} = \text{mean } F_\lambda.$$

Since the variance is $O(n)$, by Chebyshev's inequality we have that $F_\lambda(na) - F_\lambda(n\nu) \to 1$. Therefore,

$$\liminf_{n \to \infty} \frac{\log P[S_n \leq an]}{n \log n} \geq \nu - 1.$$

But $\nu$ is arbitrarily close to $a$, so the result is proved. $\quad \square$

**4. The uniform measure on $\partial B(an)$.** Let $\{\zeta_i', 1 \leq i \leq n\}$ have the distribution of $\{\zeta_i, 1 \leq i \leq n\}$ conditional on $\sum_{i=1}^{n} \zeta_i = \lfloor an \rfloor$. Let $\{\zeta_i^{(\lambda)}, 1 \leq i \leq n\}$ be independent with distribution

$$dF_{\lambda,i}(x) = \frac{1}{\phi_i(\lambda)} e^{-\lambda x} \, dF_i(x),$$

where $F_i$, $\phi_i$ are respectively the distribution function and the Laplace transform of $\zeta_i$, and $\lambda$ is the optimal parameter of the previous section, $e^{-\lambda} = b/n$. It is easy to see that $\zeta_i^{(\lambda)}$ is another Bernoulli random variable with

$$P[\zeta_i^{(\lambda)} = 1] = P[\zeta_i = 1]e^{-\lambda} \frac{1}{\phi_i(\lambda)} = \frac{1}{1 + n/(b(i-1))} := \beta_i.$$

We are now ready to prove:

THEOREM 3. *Let $|\mathcal{C}_1|$ be the length of the cycle that contains 1. Under $\mu$, the uniform distribution on $\partial B(I, an)$,*

$$|\mathcal{C}_1| \Rightarrow \mathbf{G},$$

*where $\mathbf{G}$ is a geometric r.v. with $P(\mathbf{G} > k) = (b/(1+b))^k$ and $b$ satisfies $\log(1+b)/b = 1 - a$.*

*Sketch of the proof.* The first part of demonstrating this is to recall what Arratia, Barbour and Tavaré [1] call the Feller coupling. Start with vertex 1 and choose $\sigma(1)$ uniformly from the $n$ possible choices. If this is 1, then take vertex 2 and choose $\sigma(2)$ uniformly from the $n-1$ remaining possible choices. If $\sigma(1) \neq 1$, then choose $\sigma(\sigma(1))$ uniformly from the $n-1$ remaining choices, and so on, until the final vertex where there is only one possible choice. Although the construction is much different from the Chinese restaurant process, the reader should note that if $\xi_i$ is defined by $\xi_i = 1$ if a cycle is not completed at the $i$th stage and 0 otherwise, then $\{\xi_i : 1 \leq i \leq n\}$ and $\{\zeta_i : 1 \leq i \leq n\}$ have the same distribution.

From the last observation it follows that $N = \inf\{i : \xi_i = 1\}$ has the same distribution as the length of the cycle containing 1. We can now conclude the proof of the theorem, using the large deviation calculation of the volume, and an argument called the Gibbs conditioning principle (see [4]). This principle asserts that the distribution of the $\zeta_i$ conditional on $\sum_{i=1}^{n} \zeta_i = an$ should be asymptotically independent and their law given by that which minimizes the entropy, that is, the random variables $\zeta_i^{(\lambda)}$ with distribution

$$(4.1) \qquad \frac{1}{\phi_i(\lambda)} e^{-\lambda x} \, dF_i(x)$$

where $F_i$, $\phi_i$ are respectively the d.f. and the Laplace transform of $\zeta_i$, and $\lambda$ is the parameter that optimizes (3.1), that is, $e^{-\lambda} = b/n$.

PROOF OF THEOREM 3. We will first need a lemma.

LEMMA 3. *For any $n \geq 1$ and for every $\lambda > 0$,*

$$(\zeta_1', \ldots, \zeta_n') \stackrel{d}{=} (\zeta_1^{(\lambda)}, \ldots, \zeta_n^{(\lambda)}) \qquad given \sum_{i=1}^{n} \zeta_i^{(\lambda)} = \lfloor an \rfloor.$$

PROOF. Let $f_1, \ldots, f_n$ be bounded nonnegative Borel functions:

$$E\left(\prod_{i=1}^{n} f_i(\zeta_i')\right) = E\left(f_1(\zeta_1) \cdots f_n(\zeta_n) \Big| \sum_{i=1}^{n} \zeta_i = \lfloor an \rfloor\right)$$

$$= E\left(f_1(\zeta_1) \cdots f_n(\zeta_n); \sum_{i=1}^{n} \zeta_i = \lfloor an \rfloor\right) P\left(\sum_{i=1}^{n} \zeta_i = \lfloor an \rfloor\right)^{-1}.$$

On the other hand,

$$E\left(f_1(\zeta_1^{(\lambda)}) \cdots f_n(\zeta_n^{(\lambda)}); \sum_{i=1}^{n} \zeta_i^{(\lambda)} = \lfloor an \rfloor\right)$$

$$= \int_{R^n} f_1(x_1) \cdots f_n(x_n) \mathbf{1}_{\{\sum x_i = \lfloor an \rfloor\}} \prod_{i=1}^{n} dF_{i,\lambda}(x)$$

$$= \int_{R^n} \prod_{i=1}^{n} f_i(x_i) \mathbf{1}_{\{\sum_i x_i = \lfloor an \rfloor\}} \prod_{i=1}^{n} \frac{e^{-\lambda x}}{\phi_i(\lambda)} \, dF_i(x_i)$$

$$= \frac{e^{-\lambda \lfloor an \rfloor}}{\prod_{i=1}^{n} \phi_i(\lambda)} \int_{R^n} \prod_{i=1}^{n} f_i(x_i) \mathbf{1}_{\{\sum_i x_i = \lfloor an \rfloor\}} \prod_{i=1}^{n} dF_i(x_i)$$

$$= \frac{e^{-\lambda \lfloor an \rfloor}}{\prod_{i=1}^{n} \phi_i(\lambda)} E\left[ f_1(\zeta_1) \cdots f_n(\zeta_n); \sum_{i=1}^{n} \zeta_i = \lfloor an \rfloor \right].$$

We can now divide and multiply by the probability of the events in the two sides of this equation to obtain that for some constant $C > 0$

$$E\left( \prod_{i=1}^{n} f_i(\zeta_i^{(\lambda)}); \left| \sum_{i=1}^{n} \zeta_i^{(\lambda)} = \lfloor an \rfloor \right. \right) = CE\left( \prod_{i=1}^{n} f_i(\zeta_i) \left| \sum_{i=1}^{n} \zeta_i = \lfloor an \rfloor \right. \right).$$

By taking $f_1 = \cdots = f_n = 1$ we see that $C = 1$ and the lemma is proved. $\quad\square$

We will need another lemma:

LEMMA 4. *The $\zeta_i^{(\lambda)}$ satisfy a local central limit theorem:*

$$P\left( \sum_{i=1}^{n} \zeta_i^{(\lambda)} = \lfloor an \rfloor \right) \sim Cn^{-1/2}.$$

PROOF. The proof of this local limit theorem follows very closely that of the usual i.i.d. case, which can be found in Theorem 5.2 of [5]. Let $\beta_m = P(\zeta_m^{(\lambda)} = 1)$ [i.e., $\beta_m = (1 + n/b(m-1))^{-1}$], and let $X_{m,n} = n^{-1/2}(\zeta_m^{(\lambda)} - \beta_m)$ be the rescaled Bernoulli variable. We start by noticing that $X_{m,n}$ satisfy the hypotheses of the Lindeberg–Feller theorem (Theorem 4.5 in [5]). Indeed, they are independent by definition; for all $\varepsilon > 0$, $P(|X_{m,n}| > \varepsilon) = 0$ as soon as $n^{-1/2} \le \varepsilon$, since $\zeta_m^{(\lambda)} \le 1$ and $\beta_m \le 1$ as well. Moreover,

$$\sum_{m=1}^{n} E(X_{n,m}^2) = \frac{1}{n} \sum_{m=1}^{n} \beta_m(1 - \beta_m)$$

$$\to \int_0^1 \frac{x/b}{(1 + (x/b))^2} \, dx := \sigma^2.$$

Therefore $\sum_{m=1}^{n} X_{m,n} \Rightarrow \mathcal{N}(0, \sigma)$. At this point, the proof of the local limit theorem from [5] can be reproduced exactly. Therefore

$$\sup_{x \in \mathbf{R}} \left| n^{1/2} P\left( \sum_{m=1}^{n} X_{m,n} = x \right) - n(x) \right| \to 0,$$

where $n(x) := (2\pi\sigma^2)^{-1/2} \exp(-x^2/2\sigma^2)$. Since $\sum_{m=1}^n \beta_m \to an$, and since $n(\cdot)$ is a continuous function, we can conclude the proof of the lemma by the above uniform convergence. $\quad\square$

Now, by the Feller coupling, $|\mathcal{C}_1| \overset{d}{=} \inf\{k \ge 1 : \zeta'_{n-k} = 0\}$, that is, we must reverse the time of the Chinese restaurant process. Hence by Lemmas 3 and 4:

$$
\begin{aligned}
P[\mathcal{C}_1 > k] &= P[\zeta'_n = 1, \dots, \zeta'_{n-k+1} = 1] \\
&= P\left[\zeta_n^{(\lambda)} = 1, \dots, \zeta_{n-k+1}^{(\lambda)} = 1 \,\Big|\, \sum_{i=1}^n \zeta_i^{(\lambda)} = \lfloor an \rfloor\right] \\
&= \frac{1}{P[\sum_{i=1}^n \zeta_i^{(\lambda)} = \lfloor an \rfloor]} P\left[\prod_{i=0}^{k-1} \zeta_{n-i}^{(\lambda)} = 1; \ \sum_{i=1}^{n-k-1} \zeta_i^{(\lambda)} = \lfloor an \rfloor - k\right] \\
&= \frac{1}{P[\sum_{i=1}^n \zeta_i^{(\lambda)} = \lfloor an \rfloor]} \prod_{i=0}^{k-1} P[\zeta_{n-i}^{(\lambda)} = 1] P\left[\sum_{i=1}^{n-k-1} \zeta_i^{(\lambda)} = \lfloor an \rfloor - k\right] \\
&\sim \frac{1}{1 + n/(b(n-1))} \cdots \frac{1}{1 + n/(b(n-k))} \\
&\to \frac{1}{(1 + 1/b)^k}.
\end{aligned}
$$

Hence Theorem 2 is proved. $\quad\square$

## 5. Asymptotic singularity between $\mu$ and $\nu$.

In this section we give a proof of Theorem 5 that follows in an almost straightforward way from Theorems 2 and 3: $\nu$ and $\mu$ concentrate on permutations that have a different number of fixed points. First recall the statement of the theorem:

THEOREM 5. *As $n \to \infty$, the hitting distribution $\nu$ and the uniform distribution $\mu$ on a sphere of radius $an$ are asymptotically singular:*

$$d_{\mathrm{TV}}(\mu, \nu) \to 1.$$

LEMMA 5. *The random partition of $\{1, \dots, n\}$ derived from $\nu$ is exchangeable.*

PROOF. The probability to obtain a certain partition of $\{1, \dots, n\}$ under $\nu$ only depends on the size of its blocs, which stays the same under the action of a given permutation. Hence $\nu$ yields an exchangeable partition of $\{1, \dots, n\}$. $\quad\square$

An immediate consequence is that the expected number of fixed points is $n\nu(\mathcal{C}_1 = 1) = n/(1 + b)$. Next we show that under $\nu$ the number of fixed points $N$ is close to its expected value.

LEMMA 6.
$$\operatorname{var} N = o(n^2)$$

*under $\nu$.*

PROOF. Let $x_i = \mathbf{1}_{\{\zeta_i'=0; \zeta_{i+1}'=0\}}$ be the indicator of the event that in the conditioned Chinese restaurant process, client number $i$ sits by himself. Then $N = \sum_i x_i$ and

$$\operatorname{var} N = \sum_{i=1}^n \operatorname{var} x_i + \tfrac{1}{2} \sum_{i<j} \operatorname{cov}(x_i, x_j)$$

$$\leq n + \tfrac{1}{2} \sum_{i<j} \operatorname{cov}(x_i, x_j).$$

But when $j - i > 1$, by the Gibbs asymptotic independence proved in Theorem 3, $\operatorname{cov}(x_i, x_j) \to 0$. Also, there are only $O(n)$ terms such that $j = i + 1$ and in this case $\operatorname{cov}(x_i, x_{i+1}) \leq 1$, hence the sum $\sum_{i<j} \operatorname{cov}(x_i, x_j) = o(n^2)$. □

To end the proof of Theorem 5 by Chebyshev's inequality there remains only to notice that:

LEMMA 7. *For $0 < a < 1$ and large enough $n$*

$$\nu(|\mathcal{C}_1| = 1) \neq \mu(|\mathcal{C}_1| = 1).$$

PROOF. Recall that $b$ is defined by $\log(1+b)/b = 1 - a$. For $x \in (0,1)$, let $f(x) = 1 - \log(1+x)/x$, so that $b = f^{-1}(a)$.

On the other hand, an easy consequence of Berestycki and Durrett [3] or Theorem 0 is $\mu(|\mathcal{C}_1| = 1) = e^{-u^{-1}(a)}$. [Indeed, under $\mu$, $|\mathcal{C}_1|$ is asymptotically the total progeny of a Poisson–Galton–Watson process, or *PGW* process with offspring mean $u^{-1}(a)$.]

Hence the lemma is proved if we show that

$$1/(1+b) \neq e^{-u^{-1}(a)} \quad \text{or} \quad u(x) \neq 1 - x/(e^x - 1)$$

for all $x > 0$.

We start by noticing that as $x \to 0$, $u(x) \sim x$ but $1 - x/(e^x - 1) \sim x/2$. Hence $u(x) > 1 - x/(e^x - 1)$ as $x \to 0$. The same is true as $x \to \infty$ [an easy argument shows indeed that $u(x) = 1 - e^{-x} + o(e^{-x})$]. Now those functions are both concave as we will see in a moment, hence this has to stay true on the whole open half-line $x > 0$. (Notice that we have thus proved that the hitting distribution has always fewer fixed points than the uniform distribution.) □

Lemma 8. *The function $u$ appearing in Theorem [0] is concave.*

Proof. For $c \leq 1$ this is obvious. When $c > 1$, rather than carrying explicit calculations on the second derivative of $u$, we use a theoretic argument that exploits the recent result of Schramm [12], which says that the sizes of the pieces of the giant component in the random graph have approximately a Poisson–Dirichlet distribution. Since each fragmentation decreases the distance by 1 and each coalescence increases it by 1, it is easy to see that

$$\frac{d}{dc}\mathbf{E}[d(\sigma_{cn/2}, I)|\mathcal{F}_{cn/2}] = 1 - 2P[\text{fragm.}|\mathcal{F}_{cn/2}],$$

where $\mathcal{F}.$ is the canonical filtration generated by the random walk. So we need to show that $P[\text{fragm.}]$ is an asymptotically increasing function of $c$. However, the probability of fragmenting a small cycle is asymptotically 0 (by duality and the fact that $u$ is linear in the subcritical regime), and the probability of fragmenting one of the giant cycles can be computed explicitly using the Poisson–Dirichlet structure:

$$P[\text{fragm.}] \to \mathbf{E}\sum_{i=1}^{\infty}(\theta(c)X_i)^2 = \theta(c)^2\mathbf{E}\sum_{i=1}^{\infty}X_i^2 = \tfrac{1}{2}\theta(c)^2,$$

where $\theta(c)$ is the survival probability of a $PGW(c)$ branching process and $(X_i, i \geq 1)$ follows the $PD(1)$ distribution. ($\mathbf{E}\sum X_i^2 = 1/2$ follows from [10], formula (128).) Since $\theta(c)$ is an increasing function of $c$, the lemma is proved (and thus, so is Theorem [5]). □

Remark. We have thus proved the formula

$$u(c) = c/2 - \int_0^{c/2} \theta(u)^2 \, du$$

which is perhaps a little simpler to handle than the expression in Theorem [0].

**6. Number of geodesics and Radon–Nikodym derivative.** Here we prove the following theorem, which we will then use to prove a stronger version of the singularity theorem.

Theorem 6. *Suppose $d(\sigma) = t$ and $m_1, \ldots, m_j$ are the cycle lengths of $\sigma$. The number of paths of length $t$ from $I$ to $\sigma$ is*

$$t!\prod_{i=1}^{j}\frac{m_i^{m_i-2}}{(m_i-1)!}.$$

*From this it follows that if $t = [cn/2]$ with $c < 1$, then*

$$r(\sigma) = K_{n,t} \prod_{i=1}^{j} \frac{m_i^{m_i-2}}{(m_i - 1)!},$$

*where $K_{n,t}$ is a constant that only depends on $n$ and $t$.*

*Sketch of the proof.* To see the first result, note that in order to go from $\sigma$ to $I$ in the shortest number of steps we must increase the number of cycles by 1 at each step, and to do this we must fragment a cycle at each step by transposing two of its elements. A cycle of length $m_i$ will require $m_i - 1$ fragmentations. The first step in constructing a path is to decide on how to allocate the $t$ moves between the original cycles, which can be done in $t! / \prod_{i=1}^{j}(m_i - 1)!$ ways. The next step is to count the number of ways that we can reduce a cycle of length $m_i$ in $m_i - 1$ steps, which turns out to be simple: $m_i^{m_i-1}$.

PROOF OF THEOREM 6. Given a partition of $\{1, 2, \ldots, n\}$ into groups $A_1, \ldots, A_j$ of sizes $m_i$, $1 \leq i \leq j$, the number of forests that consist of trees with vertex sets $A_1, \ldots, A_j$ is by Cayley's formula for the number of unrooted trees on $m_i$ vertices

$$\prod_{i=1}^{j} m_i^{m_i-2}.$$

Let $t = \sum_i (m_i - 1)$. A given forest can be built up in $t!$ ways so there are

$$t! \prod_{i=1}^{j} m_i^{m_i-2}$$

paths for our random graph process that end up producing a given partition. The number of permutations that correspond to a given partition is

$$\prod_{i=1}^{j} (m_i - 1)!.$$

An equal number of paths end at each permutation with cycle sizes $m_i$, $1 \leq i \leq j$, so the number of paths to a given permutation is

$$t! \prod_{i=1}^{j} \frac{m_i^{m_i-2}}{(m_i - 1)!}.$$

If $t = [cn/2]$ with $c < 1$, then the number of edge choices that end up producing no fragmentations is by Theorem 1 in [3]

$$\sim \left( \binom{n}{2} \right)^t e^{-\kappa(c)}.$$

Taking the ratio of the last two results gives Theorem 6. □

**7. The size of the support of the hitting distribution.**  In this section we prove Theorem 7, restated below.

THEOREM 7.  *Suppose $a < 1/2$. There exists a set $S_n \in \partial B(an)$ such that $\nu(S_n) \to 1$ as $n \to \infty$ and*

$$\lim_{n \to \infty} \frac{1}{n} \log \frac{|S_n|}{|\partial B(an)|} = \gamma < 0.$$

*Sketch of the proof.*  Obtaining a decay at least exponential is not very hard, even in the case $a > 1/2$. However, it is not easy to prove that this is the correct rate for the decay of $|S|/|\partial B|$, and we restrict ourselves to the case $a < 1/2$.

If $\sigma \in \partial B(an)$, then we can use Theorem 6 to find that

$$\log \nu_0(\sigma) = -an \log n + an + \sum_{k=1}^{n} a_k \log p_k + o(n),$$

where $p_k$ is the Borel distribution with parameter $c$, and $a_k$ is the number of cycles of $\sigma$ of size $k$. But by the law of large numbers, $\nu_0(a_k/n)$ should have a limit as $n \to \infty$. Hence there is a set $S$ such that $(\log \nu_0(\sigma) + an \log n)/n$ has a limit $-c_1$ whenever $\sigma \in S$. Because $\nu_0(S) \approx 1$, $|S| \approx \exp(an \log n + c_1 n)$. Moreover it is also true that $\nu(S) \to 1$. On the other hand, precise estimates on the size of $\partial B(an)$ obtained via Kolchin's representation theorem tell us that $|\partial B(an)| = \exp(an \log n + c_2 n + o(n))$. (A statement of Kolchin's theorem can be found below.) Thus, the theorem holds with $\gamma = c_1 - c_2$. To prove that $\gamma \neq 0$, we argue that the decay has to be at least exponential (a consequence of Kolchin's representation theorem).

PROOF OF THEOREM 7.  We will need precise estimates on the size of $\partial B(an)$. Because we need estimates to order higher than just $n \log n$, sticking to the large deviation approach is not good enough. Rather, we will use Kolchin's representation theorem. We would like to thank Jim Pitman for pointing out this reference to us.

Suppose we can partition $\{1, \ldots, n\}$ into a certain number of clusters, which can all have different internal states. To be more specific, suppose that each partition of $\{1, \ldots, n\}$ into $k$ clusters leads to $v_k$ possible global states of the system $\{1, \ldots, n\}$, and that we can further assign each cluster of size $j$ one of $w_j$ possible internal states. We call such a combinatorial structure a $(v, w)$-partition (of $\{1, \ldots, n\}$). Kolchin's representation theorem answers with probabilistic means to the following purely combinatorial question: how many different $(v, w)$-partitions are there? Also, what does a random, uniform, $(v, w)$-partition look like?

Before going into the details of this theorem, let us see its relevance to our problem. The number of permutations at distance $an$ from the identity is a special instance of the above Kolchin problem, where $v_k = \mathbf{1}_{\{k=(1-a)n\}}$ and $w_j = (j-1)!$. Indeed a permutation at distance $an$ is exactly a permutation having $(1-a)n$ cycles and each cluster of size $j$ can be in one of the $(j-1)!$ possible orderings of the cycle.

Here is the content of Kolchin's theorem (see [10]). Let $v(\theta) = \sum_{k=1}^{\infty} v_k \theta^k / k!$ and let $w(\xi) = \sum_{j=1}^{\infty} w_j \xi^j / j!$ be the so-called exponential generating function of the sequences $v$ and $w$. Let $K$ be an integer-valued random variable with distribution

$$P(K = k) = v_k \frac{w(\xi)^k}{k! v(w(\xi))}$$

and let $X$ be a random variable distributed according to

$$P(X = j) = \frac{w_j \xi^j}{j! w(\xi)}.$$

Here $\xi$ is any parameter. In our setting, $K = (1-a)n$, a.s. and $X$ has the so-called logarithmic distribution, $P(X = j) = b^j / j \cdot \frac{1}{-\log(1-b)}$, for some parameter $b = w(\xi)$.

THEOREM 9 (Kolchin). *The number of $(v, w)$-partitions is given by*

$$\frac{n! v(w(\xi))}{\xi^n} P\left( \sum_{i=1}^{K} X_i = n \right),$$

*where $X_i$ are i.i.d. samples of the variable $X$. Moreover, the sizes of the clusters in exchangeable random order have the same law as*

$$(X_1, \ldots, X_K) \qquad \text{given } X_1 + \cdots + X_K = n.$$

For a precise definition of exchangeable random order, and further discussion of this theorem, see [10]. It is to be noted that here $\xi$ is any parameter. By playing on this parameter so as to make the event $S_K = n$ not unlikely (e.g., of probability $\propto n^{-1/2}$ rather than exponentially small), we get that the sizes of the clusters are approximately drawn from the r.v. $X$. Note that as a consequence we get here another proof of Theorem 3. Indeed, we see that the sizes of the cycles of a uniform permutation on $\partial B$ in exchangeable random order have a logarithmic distribution (asymptotically when the parameter is chosen suitably). Hence, a size-biased pick $|\mathcal{C}_1|$ should have distribution $P(X' = j) = \text{const.} j \cdot \frac{b^j}{j} \propto b^j$, a geometric random variable. The similarity between the large deviations–statistical mechanics approach and Kolchin's theorem is striking.

Another straightforward consequence of this theorem is the precise asymptotics for the size of a ball of radius $an$. Indeed, in our setting, $v(\theta) = \theta^{(1-a)n}/[(1-a)n]!$ and $w(\xi) = -\log(1-\xi)$, hence:

$$|\partial B(an)| = \frac{n!(-\log(1-\xi))^{(1-a)n}}{\xi^n((1-a)n)!} P\left(\sum_{i=1}^{(1-a)n} X_i = n\right),$$

where $\xi$ is still any parameter. However, when $\xi$ is chosen such that $(1-a) \times E(X) = 1$, the local central limit theorem shows that $P(\sum_{i=1}^{(1-a)n} X_i = n) \sim Cn^{-1/2}$. By Stirling's formula, it is now straightforward to see that

$$|\partial B(an)| = \exp(an\log n + c_2 n + o(n)).$$

Let us now turn our attention to the hitting distribution. We will get the corresponding estimate by analyzing the Radon–Nikodym derivative $r(\sigma)$ and the law of large numbers for $\nu$, as mentioned in the sketch of the proof.

More precisely, it follows from the proof of Theorem 6 that if $\sigma \in \partial B(an)$, with cycle decomposition of size $m_1, \ldots, m_{(1-a)n}$, and $t = an$, then

$$\nu_0(\sigma) = \frac{1}{\binom{n}{2}^t e^{-\kappa(c)}} t! \prod_{i=1}^{n(1-a)} \frac{m_i^{m_i-2}}{(m_i-1)!}.$$

Let us write $a_k$ for the number of cycles of $\sigma$ of size $k$, so that $\sum_{k=1}^n a_k = n(1-a)$ and $\sum_{k=1}^n ka_k = n$. We can rewrite the above as

$$\nu_0(\sigma) = \frac{1}{\binom{n}{2}^t e^{-\kappa(c)}} t! \prod_{k=1}^n \left(\frac{1}{c}\frac{k^{k-2}}{(k-1)!}(ce^{-c})^k\right)^{a_k} c^{(1-a)n}\frac{1}{(ce^{-c})^n}.$$

When we take the logarithm, calling $q_k = \frac{1}{c}\frac{k^{k-2}}{k!}(ce^{-c})^k$ and $p_k = kq_k$,

$$\log \nu_0(\sigma) = \kappa(c) - t\log\binom{n}{2} + \log t!$$

$$+ \sum_{k=1}^n a_k \log p_k + n(1-a)\log c - n\log c + cn.$$

Recalling that $t = an$, $c = 2a$, and using Stirling's formula, we find that

$$(7.1) \qquad \log \nu_0(\sigma) = -an\log n + an + \sum_{k=1}^n a_k \log p_k + o(n).$$

We would now like to use the law of large numbers for $\nu_0$ since we know that under $\nu$,

$$\nu\left(\left\{\sigma \in \partial B : \left|\frac{a_k}{n} - q_k\right| \le \varepsilon; \ \forall\, 1 \le k \le n\right\}\right) \to 1$$

for given $\varepsilon > 0$, and where $p_k$ is the Borel distribution. But it is not directly possible to take $S = \{\sigma \in \partial B : \left|\frac{a_k}{n} - q_k\right| \leq \varepsilon; \ \forall 1 \leq k \leq n\}$ since we would obtain a bound $\varepsilon \sum_{k=1}^{\infty} \log p_k = \infty$. So we need to modify our choice: let

$$S_n = \left\{\sigma \in \partial B : \left|\frac{a_k}{n} - q_k\right| \leq (\log n)^{-5}; \ \forall 1 \leq k \leq (\log n)^2 \text{ and } a_k = 0 \text{ otherwise}\right\}.$$

There are two things we need to check on $S_n$. First we need to see that it has the property that $\nu(S_n) \to 1$, and also that it has the correct size asymptotically. The first thing is taken care of by the next lemma, while the second will follow from the fact that $\nu_0(S_n) \to 1$, itself also a consequence of the lemma below.

LEMMA 9. $\nu(S_n) \to 1$.

PROOF. It suffices to prove that $\nu(\partial B - S_n) \to 0$. By the coupling with an Erdős–Renyi random graph, there is a $\beta > 0$ such that no cycle can be greater than $\beta \log n$ with high probability under $\nu$. Hence it remains to prove that

$$\sum_{k=1}^{(\log n)^2} \nu\left(\left|\frac{a_k}{n} - q_k\right| \leq (\log n)^{-5}\right) \to 0.$$

The basic idea is to use random graph estimates. Let $G(t)$ be the result of a random graph process where edges are added in a Poissonian way at rate $\binom{n}{2}$. Then the expectation of the number of clusters of size $k$ $a'_k$ in $G(an)$, is known to be $nq_k$ asymptotically with standard deviation $O(n^{1/2})$, which is much smaller than the $n(\log n)^{-5}$ from the definition of $S_n$. So we need to show this still holds under $\nu$. Recall that we can couple the process $(G(t), t \geq 0)$ with a random walk $(\sigma_t, t \geq 0)$ where we multiply by a transposition $(i, j)$ whenever edge $(i, j)$ arrives in $G(t)$. Thus we may consider the first time $T$ that $\sigma$ is at distance $\lfloor an \rfloor$ of the identity, and obviously a realization of $\nu$ is obtained as $\sigma_T$. We consider also $T'$ the first time that $\lfloor an \rfloor$ edges have been added to $G(t)$. Then since $T'$ is the first time a Poisson process with intensity 1 exceeds the value $\lfloor an \rfloor$, $T'$ has mean $\lfloor an \rfloor$ and variance $O(n)$. Hence the number of edges added to $G$ between $\lfloor an \rfloor \wedge T'$ and $\lfloor an \rfloor \vee T'$ is $O(n^{1/2})$. On the other hand, only a bounded number $Z_{an}$ of edges are added between $T'$ and $T$, corresponding to the number of fragmentations at time $an$, and thus all in all only $O(n^{1/2})$ edges are added between $an \wedge T$ and $an \vee T$, so this may create or destroy at most $O(n^{1/2})$ clusters of size $k$ for each $1 \leq k \leq (\log n)^2$, which is much smaller than the $n(\log n)^{-5}$ from the definition of $S_n$. Hence we conclude that

$$\sum_{k=1}^{(\log n)^2} \nu\left(\left|\frac{a_k}{n} - q_k\right| \leq (\log n)^{-5}\right) \to 0.$$

□

Thus $\nu(\partial B - S_n) \to 1$. Since $\nu_0$ is obtained as an asymptotically nondegenerate conditioning of $\nu$, it follows immediately that $\nu_0(\partial B - S_n) \to 0$ as well, that is, $\nu_0(S_n) \to 1$. We now show how to use this to estimate the size of $S_n$. By (7.1), for all $\sigma \in S_n$,

$$\frac{\log \nu_0(\sigma) + an \log n}{n} \geq a + \sum_{k=1}^{(\log n)^2} (q_k + (\log n)^{-5}) \log p_k + o(1)$$

$$\geq a + \sum_{k=1}^{(\log n)^2} q_k \log p_k + o(1),$$

from which we deduce that

$$\liminf_{n \to \infty} \frac{1}{n} (\log \nu_0(\sigma) + an \log n) = a + \sum_{k=1}^{\infty} q_k \log p_k := -c_1.$$

After similar treatment for the lim sup, we get

$$\frac{1}{n} (\log \nu_0(\sigma) + an \log n) \to -c_1.$$

Since

$$\nu_0(S_n) = \sum_{\sigma \in S_n} \nu_0(\sigma) \to 1,$$

it must be that $|S_n| = \exp(an \log n + c_1 + o(n))$. Therefore

$$\lim_{n \to \infty} \frac{1}{n} \log \frac{|S_n|}{|\partial B|} = c_1 - c_2 := \gamma.$$

It now remains to show that $\gamma \neq 0$. Observe that by Kolchin's theorem, we could pursue the asymptotic expansion of $|\partial B(an)|$ and the next term would be polynomial in $n$. From the exact formula of $\nu(\sigma)$, we could also find the next term for $|S|$ and find that it is polynomial. Hence if $\gamma = 0$, then we would have $|S|/|\partial B| \sim n^{-\alpha}$ for some $\alpha \geq 0$. But, another consequence of Kolchin's theorem is that the decay has to be at least exponential: for instance, permutations in $S$ have a number of fixed points characteristic of $\nu$ and not of $\mu$. As we have seen earlier the number of fixed points under $\mu$, $n/(1+b)$, is smaller than under the hitting distribution. But since the number of fixed points under $\mu$ is given by a sum of almost independent random variables,

$$\sum_{i=1}^{n-\lfloor an \rfloor} \mathbf{1}_{\{X_i=1\}} \qquad \text{given} \qquad \sum_{i=1}^{n-\lfloor an \rfloor} X_i = n,$$

we have that

$$\mu(S) \leq P\left(\frac{1}{n} \sum_{i=1}^{n-\lfloor an \rfloor} \mathbf{1}_{\{X_i = 1\}} < \frac{1}{1+b} \,\bigg|\, \sum_{i=1}^{n-\lfloor an \rfloor} X_i = n\right)$$

$$\leq C n^{1/2} e^{-n\rho} \leq C' e^{-n\rho'}$$

by standard large deviations (here, simply Markov's inequality), and because the event on which we condition is of probability $C n^{-1/2}$. Hence the decay has to be at least exponential and $\gamma$ cannot be 0. $\square$

REMARK. The same argument shows that the hitting distribution of $\sigma$ is supported on a set at least exponentially small even in the case $a > 1/2$, but of course we do not know whether this is a precise asymptotics. If the decay is still exponential after $a = 1/2$, it seems likely that the exponential coefficient will not be smooth at $a = 1/2$. In Figure 2 we have plotted the value of this coefficient against a time-change of $a$. It would be interesting to compute exact asymptotics in the case $a > 1/2$ and make this picture complete.

REMARK. Kolchin's representation theorem could have been used already earlier for the proofs of Theorems 2 and 3. This would actually simplify the proof of both results. However, we have chosen to keep the proofs as they were, because they do not rely on a technical result such as Kolchin's theorem, which is not as well known as standard large deviation theory.

## 8. Asymptotic hyperbolicity under the uniform measure.
Here we present a proof of Theorem 8. The sketch of the proof below contains some ideas that will be used and not re-explained in the actual proof that follows.

THEOREM 10. Let $0 < a < 1$ and let $\sigma, \pi$ be two random independent points chosen uniformly from $\partial B(an)$. Then:

1. If $a < 1 - \log 2$,

$$E(\sigma|\pi)_p \leq \delta(\log n)^2$$

for some $0 < \delta = \delta(a) < \infty$. Moreover, with probability asymptotically 1, there is a geodesic between $\sigma$ and $\pi$ that comes within distance at most $\delta(\log n)^2$ of $p$.

2. If $a > 1 - \log 2$,

$$E(\sigma|\pi)_p \sim \delta n$$

for some $\delta = \delta(a) > 0$. Moreover, no geodesic can approach $p$ closer than $\delta' n$ for some $0 < \delta' < \infty$.

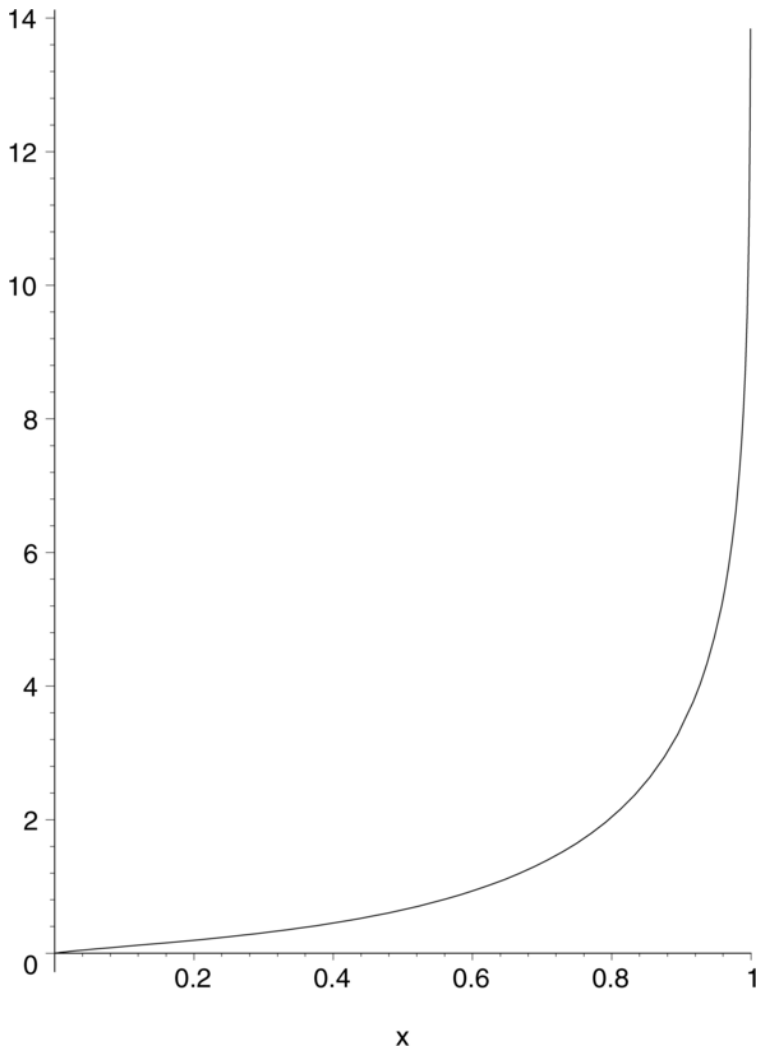

FIG. 2. *Numerical evaluation by Maple of the limiting behavior of* $-n^{-1}\log\frac{|S_n|}{|\partial B|}$. *The result is plotted as a function of* $\xi = f^{-1}(a)$, *where* $f(\xi) = 1 + \log(1-\xi)(1-\xi)/\xi$. $f^{-1}$ *is an increasing function of* $a$. *This picture has been rigorously proved only for* $a < 1/2$, *that is,* $\xi < f^{-1}(1/2) \approx 0.715331863\ldots$.

*Sketch of the proof.*   To guess what the answer is, we exploit once again the connection with the theory of random graphs. The first thing to do is to realize that because of the symmetries of the Cayley graph $G_n$ it is enough to look at $d(I, \sigma \cdot \pi)$ and see whether it is approximately $2an$ or much smaller than $2an$.

To construct our graph, we will need some notation. Let

(8.1)                    $$\pi = \tau_1 \tau_2 \cdots \tau_{an}$$

be a minimal decomposition of $\pi$ as a product of transpositions, with the following convention. If we list all cycles of $\pi$ in the order of their least element, then the transpositions $\tau_i$ are those $(x, y)$ such that $y$ comes just after $x$ in the cyclic decomposition of $\pi$ and $x$ and $y$ are in the same cycle, and we order the $an$ transpositions according to their position in this canonical decomposition. To clarify the ideas, suppose

$$\pi = (1\ 4\ 3\ 7)(2)(5\ 8)(6\ 10\ 9);$$

then we write

$$\pi = (1\ 4)(4\ 3)(3\ 7)(5\ 8)(6\ 10)(10\ 9).$$

We define the graph $\Gamma = (V, E)$ on $n(1 - a)$ vertices as follows. Let $V = \{\text{cycles of } \sigma\}$, and there is an edge between $\mathcal{C}$ and $\mathcal{C}'$ if there is $x \in \mathcal{C}$ and $y \in \mathcal{C}'$ such that $(x, y)$ is one of the $an$ transpositions in the minimal decomposition described above. Note that this graph could have self-loops and multi-edges.

A notion that we will use on several occasions is that of being a terminal point. We say $x \in \{1, \ldots, n\}$ is *terminal* if $x$ does not appear more than once in the transpositions of the above minimal decomposition. This means that, with those conventions, $x$ is situated at the "end" of the cycle of $\pi$ in which it is contained.

Here is why we are interested in the properties of the graph $\Gamma$. If we define $\sigma_0 = \sigma$ and, for $1 \leq r \leq an$, $\sigma_r = \sigma \cdot \tau_1 \cdots \tau_r$, and consider the process $(\sigma_r, 0 \leq r \leq an)$, this is a walk on $G_n$ starting at $\sigma$ and ending at $\sigma \cdot \pi$. Moreover, since at each step we are multiplying by a transposition, the cycles of $\sigma_r$ evolve according to a discrete coagulation–fragmentation chain, with cycles merging when the transposition involves elements from different cycles, and cycles splitting otherwise (as it is the case for simple random walk on $G_n$). Therefore, $\Gamma$ is the graph that results from drawing an edge between two cycles of $\sigma$ as we encounter a transposition joining those two cycles. In particular, the same argument that shows that the Erdős–Rényi random graph is an upper bound for the sizes of the cycles of simple random walk on $G_n$, will show that the cycles of $\sigma \cdot \pi$ are subcomponents of the connected components of $\Gamma$, with possibility of fragmentation whenever there is a cycle in $\Gamma$, or a self-loop or a multiple-edge. All other edges represent coalescence of cycles in the walk $(\sigma_r, 0 \leq r \leq an)$.

In particular, the property of $\Gamma$ that we will be most interested in, will be to decide whether $\Gamma$ has a giant component, meaning a component containing a positive fraction of all $n(1 - a)$ vertices. Indeed when all cycles of $\Gamma$ are small, we should expect very few cycles in $\Gamma$ and hence little fragmentation. Hence most steps of the walk $\sigma_r$ are coalescence events and the number of cycles decreases linearly; in other words, in the case that all cycles are small, $d(\sigma, \pi) \approx 2an$. On the other hand, if $\Gamma$ contains a giant cycle, then we can expect many cycles in the graph and hence many fragmentation events in the walk $(\sigma_r, 0 \le r \le an)$, which means that $d(\sigma, \pi) \ll 2an$.

Here is our strategy to see whether there is a giant component in $\Gamma$. Rather than counting the number of cycles of $\sigma$ that a component of $\Gamma$ contains, we prefer to compute the exact number of integers in $\{1, \ldots, n\}$ that it actually encloses prior to shrinking all cycles of $\sigma$ into points. Formally, this means, give weight $W(\mathcal{C}) = |\mathcal{C}|$ to any vertex $\mathcal{C}$ of $\Gamma$, and ask what is the total weight of a connected component of $\Gamma$. Let $\mathcal{C}_1(\Gamma)$ denote a size-biased pick from the connected components of $\Gamma$, that is, the total weight of the connected component of $\Gamma$ containing "1" [or, more precisely, $\mathcal{C}_1(\sigma)$].

Lemma 10 shows that $W(\mathcal{C}_1(\Gamma))$ converges in distribution to the total progeny of a branching process with offspring distribution a shifted geometric random variable. The idea is that by Theorem 3, the various cycles of $\sigma$ are asymptotically i.i.d., so that each edge in $\Gamma$ adds to the weight of $\mathcal{C}_1(\Gamma)$ a contribution which is, by Theorem 3, asymptotically a geometric random variable $\mathbf{G}$ with parameter $1/(1+b)$ where $b$ satisfies $\log(1+b)/b = 1 - a$. This seems to give an infinite progeny almost surely (since $\mathbf{G} \ge 1$ a.s.). However, to every point that we examine there is a positive probability that it is a *terminal* point. In this case, that integer does not connect to a new independent cycle of $\sigma$, and hence its offspring is 0. This kind of modified branching process is defined more precisely and analyzed in Section 8.3. The key fact is that because of the special properties of the asymptotic law of $\pi$, which involves geometric random variables, this modified branching process is in fact equal in distribution to another branching process where the offspring distribution has been shifted from $\mathbf{G}$ to $\mathbf{G} - 1$. In all that follows, we call $X$ a random variable such that

$$X \overset{d}{=} \mathbf{G} - 1.$$

Hence $\Gamma$ has a giant component if, and only if, $E(\mathbf{G}) > 2$. Since $p = 1/(1+b)$ and $\log(1+b)/b = 1 - a$,

$$P(\Gamma \text{ has a giant component}) > 0 \quad \Longleftrightarrow \quad a > 1 - \log 2.$$

PROOF OF THEOREM 8.

8.1. *Structure of the proof.* As this proof is rather long, we feel that it is appropriate to explain how the various arguments are used. In Section 8.2, we prove that $W(\mathcal{C}_1(\Gamma)) \Rightarrow \sum_{t \ge 0} Z_t$ the total progeny of a branching process with offspring distribution $X$. Then in Section 8.3, we define a modified branching process, and prove that in the case of geometric random variables this becomes another branching process with shifted offspring distribution. We then use this to prove by hand that in the subcritical case, $|\mathcal{C}_1(\sigma \cdot \pi)|$ is dominated by such a modified branching process. Since this is a subcritical branching process, we prove the exponential decay of the tail of $|\mathcal{C}_1(\sigma \cdot \pi)|$, uniformly in $n$. This enables us to show that as long as $a < 1 - \log 2$ there are very few fragmentation events in the walk $(\sigma_r, r = 0, \dots, an)$. The supercritical case is treated in Section 8.5. Since we have established branching process asymptotics, we can use the duality principle of a branching process between the subcritical phase and the supercritical phase. This shows that the number of clusters of $\Gamma$ in the supercritical regime can be computed by looking at the number of clusters of $\Gamma$ for some specific subcritical time. Since we have proved that the distance is linear in this regime, we now know how many clusters $\Gamma$ has at any subcritical time, and it follows that the number of clusters of $\Gamma$ in the supercritical regime is strictly less than what it would be if the distance was still linear. It only remains to prove that at any given time the number of extra cycles that were generated by some fragmentation (and have not been reabsorbed by other large cycles) is $O(n^{1/2})$, which is done in Section 8.6.

8.2. *Branching process asymptotics.* To start proving things, we need some more notation. Let $A_0^n = \{1\}$ and define recursively the $A_k^n$ by

$$A_{k+1}^n = \bigcup_{x \in A_k^n} \{\mathcal{C}_{\pi(x)}(\sigma)\} - \bigcup_{1 \le j \le k} A_j^n.$$

The $A_k^n$ correspond to growing the branching process generation after generation, rather than cycle after cycle. Let $(Z_t, t = 0, 1, \dots)$ be a branching process with offspring distributed as $X$. Note that by the construction of $\Gamma$, we also have that

$$\sum_{k=0}^{\infty} |A_k^n| = W(\mathcal{C}_1(\Gamma)).$$

Lemma 10. *As* $n \to \infty$,

$$(|A_0^n|, |A_1^n|, \dots) \Rightarrow (Z_0, Z_1, \dots).$$

Proof. Let us start by the convergence of $(|A_0^n|, |A_1^n|)$. If $j = 0$, $P(|A_0^n| = 1, |A_1^n| = 0) = P(\pi(1) = 1) \to 1/(1 + b) := p = P(X = 0)$. If $j \ge 1$, then

$$\begin{aligned}
P(|A_0^n| = 1, |A_1^n| = j) &= P(\pi(1) \ne 1; |\mathcal{C}_{\pi(1)}(\sigma)| = j) \\
&= P(\pi(1) \ne 1) \cdot P(|\mathcal{C}_{\pi(1)}(\sigma)| = j | \pi(1) \ne 1) \\
&\to (1 - p) \cdot (1 - p)^{j-1} p = P(X = j).
\end{aligned}$$

Indeed, conditionally on $\{\pi(1) \neq 1\}$, $\pi(1)$ is uniform on $\{2, \ldots, n\}$, so that $\mathcal{C}_{\pi(1)}(\sigma)$ is as good a size-biased pick as $\mathcal{C}_1(\sigma)$, and we can apply Theorem 3.

Now let us consider the general case of finite-dimensional distributions. Let $n_1 > 0, n_2 > 0, \ldots, n_k \geq 0$ with $\sum_i n_i \leq n$. We are trying to compute the asymptotics of

$$P(|A_0^n| = 1, |A_1^n| = n_1, \ldots, |A_k^n| = n_k).$$

To do this, we need to evaluate the probability of a collision occurring in the first $k$ stages, that is,

$$P\left(\pi(x) \in \bigcup_{1 \leq i \leq k} A_i^n - \mathcal{C}_x(\pi) \text{ for some } x \in A_j^n \text{ with } j \leq k\right).$$

We will say of an $x$ such as in the event above, that it makes a *backward* connection. Hence an $x$ makes a backward connection if $\pi(x)$ maps it to some lower level in the branching process, but $x$ is not a terminal point. Therefore backward connections (or collisions) are exactly those that may lead to a fragmentation, as explained in the sketch of the proof.

It is easy to see that $P(\text{collision in first } k \text{ stages}) = O(1/n)$. In fact, it follows from the uniformity of Lemma 11 that

$$P(|A_0^n| = n_0, \ldots, |A_k^n| = n_k; \text{ b.w. collisions in } k \text{ first stages}) \leq \sum_{j=1}^k n_j \frac{\sum_{i=0}^{j-1} n_i}{n}$$

$$\leq \frac{(\sum_{i=0}^k n_i)^2}{n}$$

(see also Lemma 12 where similar estimates are derived).

Therefore it is enough to consider

$$P(|A_0^n| = 1, |A_1^n| = n_1, \ldots, |A_k^n| = n_k| \text{ given no b.w. connection}).$$

Suppose $A_{k-1}^n = \{x_1, \ldots, x_{n_{k-1}}\}$. Conditionally on the event that there is no collision in the $k$ first stages, $\pi(x_1), \ldots, \pi(x_{n_{k-1}})$ belong to yet unexplored cycles (as long as they are not terminal). After decomposition on the number of such $x$ [call $\mathrm{T}(A_{k-1}^n)$ the number of terminal elements in the set $A_{k-1}^n$], the last probability is equal to

$$= \sum_{j=1}^n P\left(\sum_{i=1}^{n_{k-1}-j} |\mathcal{C}_{\pi(x_i)}(\sigma)| = n_k\right) P(\mathrm{T}(A_{k-1}^n - j| \text{ no b.w. conn.}))$$

$$\rightarrow \sum_{j \geq 0} P\left(\sum_{i=1}^{n_{k-1}-j} X_i = n_k\right) p^j$$

$$= P\left(\sum_{i=1}^{n_{k-1}} X_i = n_k\right)$$

by the asymptotic independence property of a finite number of size-biased cycles, and the fact that given there was no backward connection, the $x$ in level $A_{k-1}^n$ belong to different cycles of $\pi$, so that the events that they are terminal are independent asymptotically.

These are the transition probabilities of a branching process with offspring distribution $X$, so the lemma is proved. $\square$

8.3. *A modified branching process.* One way to formalize the idea that a vertex has a geometric number of children only during finitely many generations, is to use a modified branching process where each individual $x$ is endowed with a nonnegative, integer-valued random variable $T(x)$, that represents the "life-time" of its family. As long as $T(x) > 0$, $x$ will keep having children according to the original offspring progeny $L$. But when $T(x) = 0$, the individual will be declared "*terminal*" and will not be allowed to have any children.

Here is a rigorous description of this modified branching process. Let $X_{t,i}$ be a collection of i.i.d. random variables with distribution $L$, a fixed distribution on the nonnegative integers (the original progeny). Let $T_{t,i}$ be i.i.d. nonnegative integer-valued random variables, distributed according to another distribution $L'$, the lifetime. Let $Z_t$ be the size of the process at time $t$ (with discrete time). Define $Z_0 = 1$, and give the root lifetime $T_{0,0}$. Then define recursively $Z_t$ by

$$(8.2) \qquad Z_{t+1} = \sum_{i=0}^{Z_t} X_{t,i} \mathbf{1}_{\{T(x_i)>0\}},$$

where $x_1, \ldots, x_{Z_t}$ are the $Z_t$ individuals of generation $t$. If $y_1, \ldots, y_{Z_{t+1}}$ are the $Z_{t+1}$ individuals of generation $t+1$, the rule that we adopt for the value of $T(y_1), \ldots, T(y_{Z_{t+1}})$ is the following. If $T(x_i) > 0$, give all $X_{t,i}$ children of $x_i$ independent lifetimes from $T_{t,i}$, except for one of its children, say $y_j$, for which $T(y_j) := T(x_i) - 1$. Rigorously, let $N_t = \#\{i : T(x_i) > 0\}$, rewrite the $x_i$'s removing the terminal ones and call them $x'_1, \ldots, x'_{N_t}$. Let $T(y_1) = T(x'_1) - 1, \ldots, T(y_{N_t}) = T(x'_{N_t}) - 1$, and let $T(y_{N_t+1}) = T_{t+1,N_t+1}, \ldots, T(y_{Z_{t+1}}) = T_{t+1,Z_{t+1}}$.

Of course we make this definition because asymptotically, $W(\mathcal{C}_1(\Gamma))$ will be well approximated by such a system, where the offspring $L$ is the size of a cycle of $\sigma$, and where $L'$ is the size of a cycle of $\pi$. Indeed, suppose we are exploring the cluster containing 1 in the graph of the superposition of $\sigma$ and $\pi$. $T(1)$ is $|\mathcal{C}_1(\pi)|$, which corresponds to the fact that after that many iterations of $\pi$ we are back to where we started and no longer add anything new to the population of the cluster. However, after one iteration say, all vertices in the first generation, other than $\pi(1)$ itself, belong to different cycles of $\pi$ with high probability. Therefore their lifetime should be an independent random variable, distributed as $L'$.

In general, the ageing branching process $(Z_t, t = 0, 1, \ldots)$, where each individual has a "lifetime" that it transmits to one of its children, is not a Markov process with respect to its own filtration $\sigma(Z_0, Z_1, \ldots)$. Indeed the size of the generation $t + 1$ depends not only on the size of generation $t$, but also on the random variables $T_x$ where $x$ is an individual of generation $t$, so one would need to add in the filtration the values of $T(x)$ for each generation.

However, a miracle happens due to the fact that the cycles of $\pi$ have (asymptotically) a geometric distribution $\mathbf{G}$. Let $p'$ be the parameter of $\mathbf{G}$: $P(\mathbf{G} = j) = (1 - p')^{j-1} p'$. Then the distribution $L'$ of the random variables $T_{t,i}$ is again $\mathbf{G}$. For $k \geq 1$, conditionally on $T > k$, $T - k$ is distributed as $\mathbf{G}$. This fact, called "lack of memory," has the following amazing consequence:

PROPOSITION 1.    *When the lifetime $L'$ is a geometric random variable, $(Z_t, t = 0, \ldots)$ is a Markovian branching process with offspring distribution $L\mathbf{1}_{\{L'>0\}}$. When $L \stackrel{d}{=} L' \stackrel{d}{=} \mathbf{G}$, this distribution is $X \stackrel{d}{=} \mathbf{G} - 1$.*

PROOF.    Let $B_{t,i}$ be Bernoulli random variables with success parameter $P(B_{t,i} = 1) = p'$. Because $\mathbf{G} =_d \inf\{t \geq 0; B_{t,i} = 1\}$, the event $\{T(x_i) > 0\}$ is the same as $\{B_{t,i} = 0\}$, so (8.2) becomes

$$(8.3) \qquad\qquad Z_{t+1} = \sum_{i=0}^{Z_t} X_{t,i} \mathbf{1}_{\{B_{t,i}=0\}}.$$

This expresses the fact that for each new vertex visited, we can take the decision of closing the cycle, independently of the past. When the cycle still has some length to be explored, then the vertex has $X_{t,i}$ children. This decision affects the law of progeny at a given vertex. The new distribution of the progeny is now, by (8.3):

$$(8.4) \qquad P(X = j) = \begin{cases} p', & \text{if } j = 0, \\ (1 - p') P(X_{t,i} = j), & \text{if } j \geq 1. \end{cases}$$

Of course, for our problem, $\sigma =_d \pi$, so both $L$ and $L'$ are distributed as $\mathbf{G}$. As can be readily checked from (8.4), the distribution of $X$ is thus a shift of $\mathbf{G}$:

$$(8.5) \qquad\qquad\qquad X \stackrel{d}{=} \mathbf{G} - 1. \qquad\qquad\qquad \square$$

8.4. *Fragmentations in the subcritical case.*    First note that if $\sigma$ is a uniform permutation on $\partial B(an)$, if we visit all points in $\{1, \ldots, n\}$ according to their order of appearance in the canonical cyclic decomposition of $\sigma$, and call this process $V_t (0 \leq t \leq n)$, then the successive points are in some sense uniformly chosen from what remains to be found, at least as long as we do not have to start a new cycle. More precisely:

LEMMA 11. *Given $(V_0, \ldots, V_t)$ and given that $V_t$ is not a terminal point,*

$$\sigma(V_t) \text{ is uniform on } \{1, \ldots, n\} - \{V_0, \ldots, V_t\}.$$

The proof of this lemma follows directly from the Feller coupling presentation of a uniform permutation on $G_n$, which also has (obviously) this property. Conditioning on the number of cycles does not change how the cycles are filled in.

LEMMA 12. *Suppose the branching process is subcritical, that is, $p > 1/2$ or (equivalently) $a < 1 - \log 2$. Then the number of fragmentations in the walk $(\sigma_r, r = 0, \ldots, n)$ is $o(n)$.*

Basically, all cycles are fairly small, so by improving our estimates on the number of collisions, we should get an $O(1)$ bound, just like in the Erdős–Renyi case. Technicality arises due to the fact that the cycles are conditioned independent random variables, and not just independent. Here is a rigorous proof.

PROOF OF LEMMA 12. We prove things in two steps.

First, we prove a uniform bound for the size of a cluster: we show that if we write $\pi = \prod_{i=1}^{an} \tau_i$, and denote $\pi_r := \prod_{i=1}^{r} \tau_i$,

$$(8.6) \qquad P(|\mathcal{C}_1(\sigma \cdot \pi_r)| > u \text{ for some } r \leq an) \leq Cn \exp(-\alpha u),$$

where $C$ and $\alpha$ are constants independent of $n$, and $u$ is any number.

Once this exponential control is proved, we can bound the number of times that one of the $\tau_r$'s will yield a fragmentation. Indeed, recall that to obtain $\sigma \cdot \pi$ we can perform successively the $\tau_r$'s on $\sigma$, and each one yields a coagulation or a fragmentation. We hence view this as a process indexed by $1 \leq r \leq an$. In the course of this process, at all times, by (8.6) applied to $u = (\log n)^2$, *no cluster* is larger than $(\log n)^2$ with overwhelming probability, so that by Lemma 11:

$$P(\tau_{r+1} \text{ yields a fragmentation}) \leq 2(\log n)^2/n.$$

There are (exactly) $an$ transpositions to perform, hence:

$$E(\#\text{frag.}) \leq 2a(\log n)^2.$$

This is already largely enough to prove Lemma 12.

We will now prove that (8.6) holds, since this is the only thing that remains to be proved. Although we have seen that in the limit each cluster is a subcritical branching process (for which such an exponential tail of the total progeny holds), when $n$ is finite there is no real branching process available to dominate $\mathcal{C}_1(\sigma \cdot \pi_r)$, essentially because the sizes of the cycles are not i.i.d. random variables. However, they are conditionally independent (cf. Theorem 9, or Theorem 3), and we will use this fact to construct a real branching process that dominates $\mathcal{C}_1(\sigma \cdot \pi_r)$, when conditioned on some mild event. This conditioning accounts for the extra factor $n$ in (8.6).

Here is how we proceed. By Kolchin's representation theorem (Theorem 9), there are random variables $(X_1, \ldots, X_{n(1-a)})$ such that the joint law of the sizes of the cycles of $\sigma$ is $(X_1, \ldots, X_{n(1-a)})$ given $\sum_i X_i = n$ (we will call $A_n$ the event that $\sum_i X_i = n$). The $X_i$'s constitute a "pool" of possible cycle sizes. Similarly, there are random variables $(Y_1, \ldots, Y_{n(1-a)})$ such that the joint law of the sizes of the cycles of $\pi$ is $(Y_1, \ldots, Y_{n(1-a)})$ given $\sum_i Y_i = n$ (let $B_n$ be the event that $\sum_i Y_i = n$).

We give an upper bound of $\mathcal{C}_1(\sigma \cdot \pi)$ in terms of the modified branching processes of Section 8.3, that uses only the $X_i$'s and the $Y_i$'s. Start with vertex 1 and choose a size-biased pick $X_1'$ of the $X_i$'s (the cycle containing 1). Put $T(1) = Y_1'$, a size-biased pick of the $Y_i$'s. Next, given $X_1' = k$, put $T(2) = Y_2', \ldots, T(k) = Y_k'$. All vertices with positive lifetime $T$ have a number of children given by a size-biased pick of the remaining $X_i$'s. They transmit their lifetime $-1$ to one of their children and the rest have lifetimes given by size-biased picks from the remaining $Y_i$'s. Then repeat the procedure until we cannot go any further (i.e., until all vertices at a given generation have lifetime $T = 0$, or until all $X_i$'s and $Y_i$'s have been picked). Call $Z'$ the total population obtained at the end of this construction.

We claim that $Z'$ dominates all stages of $\mathcal{C}_1(\sigma \cdot \pi_r)$, because $Z'$ gives the cycles of $\sigma$ coagulated by those of $\pi$, without taking any account of eventual fragmentations. In particular, in $Z'$, as long as a vertex $x$ is not *terminal* ($T(x) > 0$), the children of $x$ will be part of the population of $Z'$. Of course in the event of a collision or a backward connection, $Z$ does not contain any additional children, so that $Z < Z'$. Therefore

$$P(Z > u) \leq P(Z' > u | A_n \text{ and } B_n)$$
$$\leq P(A_n)^{-1} P(B_n)^{-1} P(Z' > u; A_n; B_n)$$
$$\leq C n P(Z' > u).$$

Indeed, by the local central limit theorem (see [5]),

$$P(A_n) = P(B_n) = P\left(\sum_{i=1}^{n(1-a)} X_i = n\right) \sim C n^{-1/2}.$$

To complete the proof, it remains to notice that size-biasing the logarithmic distribution of Kolchin's theorem gives a geometric random variable. Therefore, by arguments already developed in the sketch of the proof, $Z'$ is the total population of a branching process with offspring distribution $X$ of (8.5), and starting with a geometric number of individuals $\mathbf{G}$. Because $p > 1/2$, this branching process is subcritical. In this case, classical estimates [2, 5] show the exponential tail

$$P(Z' > u) \leq C \exp(-\alpha u).$$

This concludes the proof of (8.6), and also that of Lemma 12.  □

REMARK.   It is possible to avoid the use of Kolchin's representation theorem in the above proof. Indeed, by Theorem 3, a size-biased pick of the cycles has, after unconditioning on some event of probability $\propto n^{-1/2}$, a distribution which is given by the lengths of sequences of 1 in the Bernoulli trials $\zeta_i^{(\lambda)}$. However, since $\beta_i = P(\zeta_i^{(\lambda)} = 1) \leq b/(1+b)$, it follows that the distribution of a size-biased cycle is thus (after unconditioning) stochastically dominated by the geometric random variable $\mathbf{G}$.

8.5. *Mean in the supercritical regime and duality.*  Although this may seem a little surprising at first, we use the result from the subcritical case to get that for the supercritical case. The idea is to use the duality of branching processes.

A crucial remark is that the number of cycles of $\sigma \cdot \pi$ is given by $\sum_{x=1}^{n} 1/|\mathcal{C}_x(\sigma \cdot \pi)|$, hence by exchangeability

$$\frac{1}{n} E(\#\text{clusters of } \Gamma) = E\left(\frac{1}{|\mathcal{C}_1(\Gamma)|}\right) \to E\left(\frac{1}{T-1}\Big| T > 1\right),$$

where $T$ is the total progeny of a branching process with offspring distributed as $X$ and started with one individual. Indeed, let us not forget that the first generation $A_0$ of the branching process is itself a geometric random variable, so we can add an imaginary root and then subtract it (thus $T - 1$). Introducing an extra vertex for the root allows us to make use of the duality principle of branching processes [2, 6].

The duality principle states that a supercritical branching process, conditioned on extinction, is another branching process, subcritical, whose offspring distribution is given through its generating function. If $\phi(s) = E(s^X)$ is the generating function of $X$ and $\alpha$ is the extinction probability $\alpha = P(T < \infty)$, then the conditioned process has offspring distribution characterized by

$$\phi'(s) = \phi(s\alpha)/\alpha.$$

Here, $P(X = j) = (1-p)^j p$, so

$$\phi(s) = \frac{p}{1 - s(1-p)}.$$

Therefore

$$\phi'(s) = \frac{p/\alpha}{1 - s\alpha(1-p)}.$$

The fixed point equation for $\alpha$ yields that

$$(8.7) \qquad \frac{p}{\alpha} + \alpha(1-p) = 1$$

so that $\phi'$ is the Laplace transform of another shifted geometric random variable $X'$, with parameter $p' = p/\alpha$. Let $T'$ be the total progeny of a branching process with offspring $X'$ and started with one individual.

Let us now relate the supercritical and subcritical regimes. By duality,

$$\begin{aligned} E\Big(\frac{1}{T-1}\Big|T>1\Big) &= \frac{P(T<\infty)}{P(T>1)}E\Big(\frac{1}{T-1}\mathbf{1}_{\{T>1\}}\Big|T<\infty\Big) \\ &= \frac{P(T<\infty)}{P(T>1)}E\Big(\frac{1}{T'-1}\mathbf{1}_{\{T'>1\}}\Big) \\ &= \frac{P(T<\infty)}{P(T>1)}P(T'>1)E\Big(\frac{1}{T'-1}\Big|T'>1\Big). \end{aligned}$$

However, for the subcritical regime, we know by Lemma 12 that there are only $o(n)$ fragmentations, so the distance between $\sigma$ and $\pi$ is $2an+o(n)$ and the number of clusters is $(1-2a)n+o(n)$. Hence $E(\frac{1}{T'-1}|T'>1) = 1-2a'$, where $a'$ is the radius corresponding to the conditioned parameter $p'$. Since $p = 1/(1+b)$ and $\log(1+b)/b = 1-a$, we have that

$$a = 1 + \frac{p\log p}{1-p}.$$

On the other hand, due to the fixed point equation (8.7), the constant $P(T'>1)/P(T>1) = (1-p')/(1-p)$ simplifies into $\alpha$.

Therefore, Theorem 8 is proved when we show that

$$\alpha^2(1-2a') > 1-2a.$$

Using the fixed point equation (8.7), we find that $a' = 1 + p\log p'/((\alpha^2)(1-p))$, so that the above reduces to

$$-\alpha^2 - 2\frac{p\log p'}{1-p} > -1 - 2\frac{p\log p}{1-p}$$

or

$$\frac{2p\log\alpha}{1-p} > \alpha^2 - 1.$$

Using one more time the fixed point equation, one gets

$$2p\log\alpha > \alpha - 1.$$

Since $\log(1-x) > -x$, it is therefore enough to show

$$-2p(1-\alpha) > \alpha - 1 \text{ or } 2p < 1$$

which is precisely the condition that the branching process is supercritical.

8.6. *Fragmentations in the supercritical range.* In the previous section we have computed asymptotics for the expected number of clusters in the graph resulting from the superposition of the cycle structures of $\sigma$ and $\pi$. We now need to show that at the end of the walk $(\sigma_r, r = 0, \ldots, an)$, there are no more than $o(n)$ additional cycles that have been generated by fragmentation, compared to the number of clusters of $\Gamma$.

To do this, we use once again the dynamic point of view adopted to deal with the subcritical regime. Let $\tau_1, \ldots, \tau_{an}$ be the decomposition of $\pi$ in product of $an$ transpositions as evoked earlier, and let $\sigma_r = \sigma \cdot \tau_1 \ldots \tau_r$.

LEMMA 13. *For each $1 \leq r \leq an$, the expected number of cycles in $\sigma_r$ generated by fragmentation is $O(n^{1/2})$.*

This is similar to the Erdős–Renyi case of Berestycki and Durrett [3], Theorem 3. Lemma 13 does not claim that the number of fragmentations itself is $O(n^{1/2})$, but that the number of extra cycles generated by fragmentation is $O(n^{1/2})$. Just like in the Erdős–Renyi case, many of the cycles that are fragmented get reabsorbed by large components fairly quickly.

PROOF OF LEMMA 13. There can never be more than $n^{1/2}$ cycles of size larger than $n^{1/2}$. On the other hand, by Lemma 11, the probability that $\tau_t$ will create a fragment of size smaller than $n^{1/2}$ is at most $n^{1/2}/n = n^{-1/2}$. Therefore the expected number of such fragmentations is at most $an \cdot n^{-1/2} = O(n^{1/2})$. □

At this point, Theorem 8 is proved. □

**Acknowledgments.** I am very grateful to Laurent Saloff-Coste for introducing me to Gromov hyperbolic spaces and to Jim Pitman for pointing out the reference to Kolchin's representation theorem. I wish to thank Rick Durrett particularly for his help at several important stages of this work and for his support during my stay at Cornell University.

Department of Mathematics
University of British Columbia
1984 Mathematics Road—Room 121
Vancouver, British Columbia
Canada V6T 1Z2
E-mail: nberestycki@math.ubc.ca